\documentclass[12pt]{article}
\usepackage{amssymb,amsmath,amscd}

\input epsf.tex

\title{Motivic Donaldson-Thomas invariants: summary of results}

\author {Maxim Kontsevich, Yan Soibelman}

\begin{document}
\maketitle

\newcommand{\op}[1]{\operatorname{#1}}

\newcommand{\CC}{{\mathcal C}}

\renewcommand{\O}{{\mathcal O}}

\newcommand{\E}{{\mathcal E}}
\newcommand{\F}{{\mathcal F}}

\newcommand{\g}{{\mathfrak g}}
\newcommand{\h}{{\mathfrak h}}

\renewcommand{\k}{{\bf k}}
\newcommand{\kk}{{\overline{\bf k}}}

\newtheorem{defn}{Definition}
\newtheorem{thm}{Theorem}
\newtheorem{lmm}{Lemma}
\newtheorem{rmk}{Remark}
\newtheorem{prp}{Proposition}
\newtheorem{conj}{Conjecture}
\newtheorem{exa}{Example}
\newtheorem{cor}{Corollary}
\newtheorem{que}{Question}
\newtheorem{ack}{Acknowledgments}
\newcommand{\C}{{\bf C}}
\newcommand{\K}{{\bf k}}
\newcommand{\R}{{\bf R}}
\newcommand{\N}{{\bf N}}
\newcommand{\Z}{{\bf Z}}
\newcommand{\Q}{{\bf Q}}
\newcommand{\G}{\Gamma}
\newcommand{\A}{A_{\infty}}

\newcommand{\epi}{\twoheadrightarrow}
\newcommand{\mono}{\hookrightarrow}
\newcommand\ra{\rightarrow}
\newcommand\uhom{{\underline{\op{Hom}}}}
\renewcommand\O{{\cal O}}
\newcommand{\epp}{\varepsilon}

\newcommand{\M}{{\mathsf{M}}}
\newcommand{\Gr}{{\mathsf{Gr}}}
\newcommand{\Fl}{{\mathsf{Fl}}}
\newcommand{\univ}{{\rm{univ}}}
\newcommand{\virt}{{\mathrm{virt}}}
\newcommand{\weak}{{\mathrm{weak}}}
\newcommand{\B}{{\mathrm{B}}}
\renewcommand{\H}{{\mathcal{H}}}

{\it To the memory of I.M. Gelfand}

\tableofcontents

\section{Introduction}

This paper is a ``digest" of  our longer paper [KoSo1] where the concept of {\it motivic Donaldson-Thomas invariants} for $3$-dimensional Calabi-Yau category ($3CY$ category for short) was introduced.

From the point of view of a physicist, in [KoSo1] we offered a mathematical definition of the notion of {\it BPS invariants} (including ``refined BPS invariants", cf. [GuDi]) in a very general, model independent, framework. In loc. cit. we also proved a general wall-crossing formula, which agrees with known wall-crossing formulas for $N=2$ theories, e.g. with the one for $SU(2)$ Seiberg-Witten theory, with Denef-Moore semiprimitive wall-crossing formula (see [DeMo]) and with Cecotti-Vafa work [CeVa] on soliton counting in $2d$ sypersymmetric QFTs.

From the point of view of a mathematician, we suggested a framework for the counting of ``motivic invariants" (e.g. Euler characteristic, Poincar\'e polynomial, etc.) of $3CY$ categories (e.g. derived category of coherent sheaves on a compact or local $3$-dimensional Calabi-Yau manifold).

Since the paper [KoSo1] is very involved and contains many constructions, ramifications and connections to many different topics, its logic can be hidden ``behind the trees". In this paper we are going to recall the general philosophy of [KoSo1] and discuss few applications and open problems.

Let us summarize the  approach of [KoSo1] which leads in the end to the definition of motivic and numerical Donaldson-Thomas (DT for short) invariants. For the motivations we refer to the Introduction in [KoSo1].

1. With an ind-constructible triangulated $\A$-category ${\cal C}$ over a field $\k$ (see [KoSo1], Section 3) we associated (see [KoSo1], Section 6.1) a unital associative algebra $H(\CC)$ which is called {\it motivic Hall algebra}. Roughly speaking it is an algebra formed by motivic stack functions (see [Jo3] and [KoSo1], Section 4.2)
on the space (stack) of objects $Ob(\CC)$. The structure constants of $H(\CC)$ are defined via an appropriate count of exact triangles in the triangulated category $\CC$.

2. A choice of constructible stability condition on $\CC$ (see [KoSo1], Section 3.4) allows us to define a collection of full subcategories $\CC_V\subset \CC$ parametrized by strict sectors $V\subset \R^2$. Each subcategory $\CC_V$ is generated by extensions of semistable objects of $\CC$ with central charges in $V$ and by the zero object. More precisely, the definition of the stability condition depends on some data. Part of the data is a homomorphism $cl:K_0(\CC)\to \Gamma$, where $\Gamma$ is a free abelian group endowed with an integer-valued bilinear form. Another piece of the data is the  {\it central charge}, i.e. homomorphism $Z:\Gamma\to \C$.
A choice of $V$ allows us to define a convex cone $C(V,Z)\subset \Gamma\otimes \R$, generated by some lattice elements $\gamma$ for which $Z(\gamma)\in V$. Then  one can complete $H(\CC_V)$ with respect to the cone and obtain a completed motivic Hall algebra $\widehat{H}(\CC_V)$.

3. We define a collection of invertible elements $A_V^{Hall}\in\widehat{H}(\CC_V)$, roughly, as characteristic functions of the stack of objects of $\CC_V$ (see [KoSo1], Section 6.1). In the case when $\k$ is a finite field one has
$$A_V^{\op{Hall}}=1+\dots=\sum_{[E]\in Iso(\CC_V)}{[E]\over {\#\op{Aut}(E)}}\,\,.$$
Elements $A_V^{Hall}$ for generic $V$ satisfy the following {\it Factorization Property}:
$$A_{V}^{\op{Hall}}=A_{V_1}^{\op{Hall}}\cdot A_{V_2}^{\op{Hall}}$$
for a clockwise decomposition of the strict sector into subsectors: $V=V_1\sqcup V_2$. Strict sectors are not required to be non-degenerate, open or closed, thus an example of the decomposition is a sector $V$ dissected into two by a ray. Furthermore, the Factorization Property implies a factorization formula
$$A_V^{Hall}=\prod_{l\subset V}^{\longrightarrow}A_l^{Hall}$$
taken in the clockwise order over all rays $l$ belonging to $V$.

4. In case when $\CC$ is an ind-constructible $3CY$ category over a field of characteristic zero (see [KoSo1], Section 3.3)  we define a {\it motivic quantum torus} ${\cal R}_{\CC}$, which is an associative algebra described by the usual relations
$$\hat{e}_{\gamma_1}\hat{e}_{\gamma_2}={\mathbb L}^{{1\over{2}}\langle \gamma_1,\gamma_2\rangle}\hat{e}_{\gamma_1+\gamma_2}, \,\,\,\hat{e}_0=1\,\,.$$
Here ${\mathbb L}^{{1\over{2}}}$ is the square root of the motive ${\mathbb L}$ of the affine line $[{\bf A}^1_{\k}]$,  elements $\gamma$ belong to the lattice $\Gamma$, and we assume that the integer-valued bilinear form $\langle \gamma_1,\gamma_2\rangle$ is skew-symmetric.
The coefficient ring for the quantum torus can be any commutative ring $C$. We use two choices for $C$: either a certain localization of the ring of motivic functions on $Spec(\k)$ or its $l$-adic version (more precisely we pass to rings formed by certain equivalence classes, see [KoSo1], Sections 4.5, 6.2).

5. After a choice of the so-called orientation data (roughly, a square root of a super line bundle $sdetExt^{\bullet}(E,E)$ on the space of objects $Ob(\CC)$, see [KoSo1], Section 5)  we define {\it motivic Donaldson-Thomas invariants} as a collection of elements $A_V^{mot}\in {{\cal R}}_{\CC_V}$ which are images of $A_V^{Hall}$ under the homomorphisms
$$\Phi_V: \widehat{H}(\CC_V)\to {{\cal R}}_{\CC_V}$$
defined in Section 6.3 of [KoSo1]. Here we slightly abuse the notation denoting by ${{\cal R}}_{\CC_V}$ a completion of ${\cal R}_{\CC}$ which depends on a choice of $V$ (see Section 4.2. below).

More precisely, we proved the existence of the collection of homomorphisms $\Phi_V$ in the $l$-adic case which is sufficient for applications to numerical DT-invariants. The general motivic case depends
on a conjectural integral formula (Conjecture 4 stated in Section 4.4 of [KoSo1]). The definition of the homomorphism $\Phi_V$ relies on the theory of {\it motivic Milnor fiber} $MF(W_E)$ of the potential $W_E$ of an object $E\in Ob(\CC)$. The potential  is defined in terms of the Calabi-Yau structure. The motivic Milnor fiber satisfies the motivic Thom-Sebastiani theorem proved by Denef and Loeser (see [DenLo]). The latter is used in the proof that $\Phi_V$ is a homomorphism of rings.

Since the elements $A_V^{Hall}$ satisfy the Factorization Property, the same is true for the elements $A_V^{mot}$.

6. Having the motivic DT-invariant (which can be informally thought of as Serre polynomial associated with the mixed Hodge module of vanishing cycles on $Ob(\CC)$ of the potential $W_E$)
we can obtain numerical DT-invariants by a procedure similar to passing from the Poincar\'e polynomial to the Euler characteristic, i.e. a kind of ``quasi-classical limit". The situation is complicated by the fact that the elements $A_V^{mot}$ defined above do not have such a limit.
We conjectured in [KoSo1], Section 7.3. that there exists a specialization of the {\it automorphism} $x\mapsto A_V^{mot}x(A_V^{mot})^{-1}$  of the (completed) motivic  quantum torus ${{\cal R}}_{\CC_V}$ at
${\mathbb L}^{{1\over{2}}}=-1$. It gives a formal Poisson automorphism $A_V$ of the Poisson torus, whose algebra of functions is generated by the elements $e_{\gamma}$ which are quasi-classical limits of $\hat{e}_{\gamma}$.

7. We define birational Poisson automorphisms of the Poisson torus associated with $\Gamma$ such as follows:
$$ T_{\gamma}: e_{\mu}\mapsto (1-e_{\gamma})^{\langle \gamma, \mu\rangle}e_{\mu},\,\,\,\,\, \gamma,\mu \in \Gamma\,\,.$$
Then each $A_V$ belongs to the pronilpotent group of formal Poisson automorphisms of the above torus and is given by the formula:
$$A_V=\prod_{\gamma\in C(V,Z)\cap \Gamma}^{\longrightarrow} T_\gamma^{\,\Omega(\gamma)}$$
where $\Omega(\gamma)$ are rational numbers, $Z$ is the central charge, and the product is taken in the clockwise direction with respect to the arguments of $Z(\gamma)$ (more precisely one defines first a factor corresponding to a single ray $l=\R_{>0}Z(\gamma)$, and then combines these factors into the clockwise product, see the details in the loc. cit). These are {\it numerical DT-invariants} of the category $\CC$ (BPS invariants  in the language of physics). We conjectured in Section 7.1 of [KoSo1] that $\Omega(\gamma)$ are integers for a generic stability condition, and we proved the integrality in [KoSo3] for quivers with potential. Also, conjecturally they do not depend on the orientation data.

The above definition of $\Omega(\gamma)$  agrees with the Behrend local invariant
$$(-1)^{\dim X}(1-\chi(MF_E(W_E)))\,\,,$$
(see [Be1]) where in the RHS we compute the Euler characteristic of the Milnor fiber of the potential $W_E$ (informally, the set of objects of $\CC$ belongs to the set of critical points of $W_E$ considered as a function on a larger space $X$). More precisely, the invariants agree for the category consisting of Schur objects. For that reason our motivic DT-invariants can be thought of as a quantization of DT-invariants introduced by Behrend in [Be1].

Summarizing, the definition of motivic DT-invariants (refined BPS invariants in the language of physics) of the category $\CC$ depends a priori on two conjectures: existence of orientation data and  the the integral identity conjecture from Section 4.4 of [KoSo1]. The definition of the numerical DT-invariants $\Omega(\gamma)$ depends on the conjecture about existence of the quasi-classical limit of motivic DT-invariants (see paragraph 6 above). It does not depend on the full version of the integral identity, since we can use the $l$-adic version of the theory (including e.g. the $l$-adic version of the quantum torus), for which the integral identity was proved in Section 4.4 of [KoSo1]. Also, we expect that numerical DT-invariants do not depend on the orientation data. As we already pointed out, our motivic DT-invariants can be informally thought of as Hilbert series of the equivariant cohomology of the motivic perverse sheaf of vanishing cycles of the potential on the stack of objects of our $3CY$ category. In some cases (e.g. for the Calabi-Yau category associated with quiver with potential) we can give more precise meaning to the above phrase. More on that can be found in our paper [KoSo3], where we offer an alternative (and more elementary) point of view on motivic DT-invariants. Although in [KoSo3] we consider only $3CY$ categories associated with representations of a quiver with potential, it seems that the approach we propose there is more general. More precisely, we expect that an ind-constructible locally regular $3CY$ category (satisfying some additional conditions)  can be locally described by a quiver with potential. Then the Hilbert series from [KoSo3] (see also Section 7.3 below)  is the motivic DT-invariant.

About the structure of the paper. It consists of two parts.
In the first (larger) part we review the constructions from [KoSo1] which lead to motivic and numerical DT-invariants.
Most of the material is borrowed from [KoSo1].  In the second part we discuss several speculations and conjectures from [KoSo1]. Some of them deserve a special project. In particular, as our original motivation for [KoSo1] was our paper [KoSo4], we stress here again that a good way to encode our DT-invariants is as gluing data for a complex integrable system (or its quantization) ``a'la [KoSo4]".
The DT-invariants are used for the definition of the gluing automorphisms. Wall-crossing formulas ensure that the gluing is well-defined. Maybe this  integrable system (or rather its quaternion K\"ahler cousin, see [KoSo5]) is the main  structure underlying the mathematical theory of BPS states.

We should say that we have not attempted to mention here contributions of all researchers working in the area of Donaldson-Thomas invariants and their generalizations. As a result our list of references contains very few papers. We do not discuss here the motivations and the history of the subject, just mentioning that  [KoSo1] was mostly motivated by our previous work [KoSo4]. At the same time, there are several papers on DT-invariants written by other authors which are logically related to [KoSo1]. We  would like to mention just few of them for the sake of the reader.

1. The approach to DT-invariants via critical points of a function is due to Behrend (see [Be1], [BeFa]). Our paper [KoSo1] can be thought of as a ``quantization" of [Be1], with the quantization parameter being the motive of affine line.

2. Joyce was probably first who approached ``DT-type" invariants in abstract categorical setting. He developed  the technique of motivic stack functions and understood the relevance of motives to the counting problem (see [Jo1-Jo4]). The main limitation of his approach was due to the fact that he did not use the notion of potential and worked with abelian rather than triangulated categories (for many  applications, especially to physics, one needs triangulated categories). His recent theory of ``generalized DT-invariants" developed jointly with Song (see [JoSo1,2]) fixes some of these gaps and fits well with the general philosophy of [KoSo1] (in fact the papers [JoSo1,2] use some parts of [KoSo1]). Those papers as well as the papers by Pandharipande and Thomas (see [PT1,2]) deal with concrete examples of categories (e.g. the category of coherent sheaves). The authors construct numerical invariants via Behrend approach. It is difficult to prove that they are in fact invariants of triangulated categories (which is manifest in [KoSo1]).

3. The general concept of BPS state when encoded mathematically uses the notion of semistable object in a triangulated category. Hence the theory of stability conditions of Bridgeland (see [Br1]) is a necessary tool. Motivated by [Jo1-4] Bridgeland developed (jointly with Toledano Laredo) an approach
to Joyce invariants based on the theory of Stokes factors for irregular connections on ${\bf P}^1$ (see [BrTL]). As we explained in Section 2.7 of [KoSo1] this is a special case of our general theory of stability data on graded Lie algebras. This part of the story is independent on the categorical framework, and we will only briefly discuss it below in Section 5 (for  more details see [KoSo1], Section 2).

{\it Acknowledgements.} Second author was partially supported by NSF FRG grant DMS-0854989. He also thanks to IHES for excellent research conditions.

\section{Ind-constructible $3CY$-categories}

For an introduction to the geometric language of $\A$-categories and axiomatics of triangulated $\A$-categories we refer the reader to [KoSo2] and the Appendix of [So2]. For a more algebraic exposition of $\A$-categories the reader can look at [Ke1]. In any case, we will assume the familiarity with some basic notions. For instance the reader should know that for an $\A$-category there are infinitely many higher composition maps $m_n, n\ge 1$, such that $m_1: Hom(E_1,E_2)\to Hom(E_1,E_2)[1]$ is a differential and the composition map $m_2: Hom(E_1,E_2)\otimes Hom(E_2,E_3)\to Hom(E_1,E_3)$ is associative only on the cohomology groups with respect to $m_1$. Differential-graded categories are examples of $\A$-categories with all $m_n=0, n\ge 3$.

For simplicity we will assume that all $\A$-categories are defined over the ground field $\k$ of characteristic zero (although this condition can be relaxed for some results, e.g. for the definition of motivic Hall algebra).

\subsection{Ind-constructible triangulated $\A$-categories, Calabi-Yau categories, stability conditions}

Recall that a  Calabi-Yau category of dimension $d$ is a weakly unital $\k$-linear triangulated $\A$-category
such that for any two objects $E,F$ the $\Z$-graded vector space
$\op{Hom}^\bullet(E,F)=\oplus_{n\in \Z}\op{Hom}^n(E,F)$ is finite-dimensional
and moreover:

1) We are given a  non-degenerate pairing
$$(\bullet, \bullet): \op{Hom}^{\bullet}(E,F)\otimes \op{Hom}^{\bullet}(F,E)\to \k[-d]\,\,,$$
which is symmetric with respect to interchange of objects $E$ and $F$;

2) For any $N\geqslant 2$ and a sequence of objects $E_1,E_2,\dots,E_N$ we are given
a polylinear $\Z/N\Z$-invariant map
$$W_N: \otimes_{1\le i\le N}\left(\op{Hom}^{\bullet}(E_i,E_{i+1})[1]\right)\to \k[3-d]\,\,.$$
Here $[1]$ means the shift in the category
of $\Z$-graded vector spaces, and we set $E_{N+1}=E_1$;

Explicitly, the maps $W_N$ are defined in terms of the higher compositions $m_n$ and the pairing by the formula
$$W_N(a_1,\dots,a_{N})=(m_{N-1}(a_1,\dots,a_{N-1}),a_N)\,\,.$$

The collection $(W_N)_{N\geqslant 2}$ is called the {\it potential of $\CC$}.
If $d=3$ then for any object $E\in Ob(\CC)$ we define a formal series $W_{E}^{tot}$ at $0\in \op{Hom}^{\bullet}(E,E)[1]$
by the formula:
$$W_E^{tot}(\alpha)=\sum_{n\geqslant 2}{W_n(\alpha,\dots, \alpha)\over{n}}\,\,.$$
We call $W_E^{tot}$ the {\it total (or full) potential of the object $E$}. We call the {\it potential of $E$}
the restriction of $W_E^{tot}$ to the subspace $\op{Hom}^1(E,E)$. We will denote it by $W_E$.
One can prove (see [KoSo1], Proposition 7) that the potential $W_E$ admits a decomposition
$$W_E=W_E^{min}\oplus Q_E\oplus N_E\,\,,$$
where $W_E^{min}$ is the potential of the minimal model $\CC^{min}$ (i.e. it is a formal series on $\op{Ext}^1(E,E)$),
the quadratic form $Q_E$ is defined on the vector space $\op{Hom}^1(E,E)/\op{Ker}(m_1:\op{Hom}^1(E,E)\to \op{Hom}^2(E,E))$ by the formula
$Q_E(\alpha,\alpha)={m_2(\alpha,\alpha)\over{2}}$, and $N_E$ is the zero function on the image of the map $m_1:\op{Hom}^0(E,E)\to \op{Hom}^1(E,E)$.
In the above splitting formula we use the notation $(f\oplus g)(x,y)=f(x)+g(y)$ for the direct sum of formal functions $f$ and $g$.

In order to apply the techniques of the theory of motivic integration we consider a class of $\A$-categories which we call ind-constructible. Recall the following definition.

\begin{defn} Let $S$ be a variety over $\k$. A subset $X\subset S(\kk)$ is called constructible over $\k$
if it belongs to the Boolean algebra generated by $\kk$-points of open (equivalently closed) subschemes of $S$.

\end{defn}

Ind-constructible sets are colimits of constructibles ones (equivalently, countable unions of non-intersecting constructible sets). They naturally form a symmetric monoidal category. The notion of ind-constructible $\A$-category developed in [KoSo1], Section 3.1 as well as its Calabi-Yau version developed in loc.cit. Section 3.3, roughly speaking, means that objects of such a category form an ind-constructible set, morphisms form constructible bundles, higher composition maps are morphisms of tensor products of such bundles, etc.
We impose the condition that near each object all the structures are in fact regular i.e. described by schemes, regular maps, etc. In particular, the potential $W_E$ can be thought of as a ``partially formal" function: it is regular along $Ob(\CC)$ near $E\in Ob(\CC)$ and formal in the direction
$\alpha\in \op{Hom}^1(E,E)$. We skip here many delicate details, referring the reader to  Section 3 of [KoSo1].

\begin{rmk} In Section 3.2 of [KoSo1] we explain how to associate to an ind-constructible category $\CC$ over $\k$ an ordinary $\k$-linear $\A$-category $\CC(\k)$. This construction is not quite trivial, since the isomorphism classes of objects in $\CC$ in the constructible setting depend only on the naturally defined category $\CC(\kk)$ over the algebraic closure $\kk$ of $\k$. The correct definition involves descent data associated with all Galois extensions of $\k$. Because of that we impose the restriction that the ground field $\k$ is perfect.

\end{rmk}

The class of ind-constructible $\A$-categories is large due to Bondal-van den Bergh theorem (see [BonVdB]), which makes many geometrically defined categories (e.g. $D^b(X)$ for a smooth projective scheme $X$) to be equivalent to the category $A$-mod of $\A$-modules over an $\A$-algebra with finite-dimensional cohomology. In this ``derived non-commutative geometry" sense all schemes become affine.  Thus one can describe  objects of a category by an increasing number of solutions to polynomial equations solved in the free modules $A^n$. This implies that the objects form and ind-constructible set. Similarly one can treat morphisms.

A modification of Bridgeland theory of stability conditions to the case of ind-constructible triangulated categories was suggested in Section 3.4 of [KoSo1]. Let us briefly recall it here.

Let $\CC$ be an ind-constructible weakly unital $\A$-category over a field $\k$ of arbitrary characteristic.
Let $\op{cl}: Ob(\CC)\to \Gamma\simeq \Z^n$ be a map of ind-constructible sets (where $\Gamma$ is considered as a countable set of points) such that the induced map $ Ob(\CC)(\overline{\k})\to \Gamma$ factorizes through a group
homomorphism $\op{cl}_{\overline{\k}}: K_0(\CC(\overline{\k}))\to \Gamma$. It is easy to see that for any field extension
$\k^{\prime}\supset \k$ we obtain a homomorphism $\op{cl}_{\overline{\k}^{\prime}}: K_0(\CC(\overline{\k}^{\prime}))\to \Gamma$.

\begin{defn} A constructible stability structure on $(\CC,\op{cl})$ is given by the following data:

\begin{itemize}
\item{ an ind-constructible subset
$$\CC^{ss}\subset Ob(\CC)$$  consisting of objects called semistable, and satisfying the condition that with each object it contains all isomorphic objects,}

\item{an additive map $Z:\Gamma\to \C$ called central charge, such that $Z(E):=Z(\op{cl}(E))\ne 0$ if $E\in \CC^{ss}$,}

\item{a choice of a branch of the logarithm $\op{Log}Z(E)\in \C$ for any $E\in \CC^{ss}$ which is constructible as a function of $E$.}

\end{itemize}
\end{defn}

These data are required to satisfy the usual Bridgeland axioms (see [Br1]). In particular, for any object we have its Harder-Narasimhan filtration. We impose also the following new axioms:

\begin{itemize}
\item{the set of $E\in \CC^{ss}(\overline{\k})\subset Ob(\CC)(\overline{\k})$ with  fixed $\op{cl}(E)\in \Gamma\setminus \{0\}$ and fixed $\op{Log}Z(E)$ is a constructible set.}

\item{(Support Property) Pick a norm $\parallel\cdot \parallel $ on $\Gamma \otimes
\R$. Then there exists $C>0 $ such that
    for all $E\in \CC^{ss}$ one has
   $ \parallel E \parallel\le C|Z(cl(E))|$.}

\end{itemize}

Slightly ahead of our exposition we remark here that
the Support Property is equivalent to the one for stability data on graded Lie algebras, which we are going to use below in Section 5.1 (see there for the details):

{\it There exists a  quadratic form $Q$ on $\Gamma_{\R}:=\Gamma\otimes \R$ such that

1) $Q_{|\op{Ker}\, Z}<0$;

2) $\op{Supp} a\subset \{\gamma\in \Gamma\setminus \{0\}|\,\,Q(\gamma)\geqslant 0\}$,

where we use the same notation $Z$ for the natural extension of $Z$ to $\Gamma_{\R}$.}
We are going to explain in Section 5 that the  notion of stability data $a=(a(\gamma))_{\gamma \in \Gamma}$ on a graded Lie algebra leads naturally to symplectomorphisms generalizing $A_V$.

\subsection{Motivic stack functions}

The theory of motivic functions introduced by the first author was developed by many people, most notably, by Denef and Loeser. Its equivariant version (theory of motivic stack functions) was developed independently by Joyce (see [Jo3]). We recall below a version based on Sections 4.1, 4.2 of [KoSo1].

The abelian group $Mot(X)$ of motivic functions is the group generated by symbols $[\pi:S\to X]:=[S\to X]$ where $\pi$ is a morphism of constructible sets, subject to the relations $$[(S_1\sqcup S_2)\to X]=[S_1\to X]+[S_2\to X]\,\,.$$
For any constructible morphism $f:X\to Y$ we have two homomorphisms of groups:

1) $f_{!}: Mot(X)\to Mot(Y)$, defined by
$[\pi:S\to X]\mapsto [f\circ \pi:S\to X]$;

2) $f^{\ast}: Mot(Y)\to Mot(X)$, defined by
$[S^{\prime}\to Y]\mapsto [S^{\prime}\times_YX\to X]$.

Moreover, $Mot(X)$ is a commutative ring via the fiber product operation.
We denote by ${\mathbb L}\in Mot(Spec(\k))$ the element $[{\bf A}^1_{\k}]:=
[{\bf A}^1_\k \to  Spec(\k)]$. It is customary to add its formal inverse
${\mathbb L}^{-1}$ to the ring $Mot(Spec(\k))$ (or more generally to the ring $Mot(X)$ which is a $Mot(Spec(\k))$-algebra).

Let us recall several ``realizations" of motivic functions. For each realization the theory of motivic DT-invariants developed in [KoSo1] has the corresponding version (except, maybe, the case (iv) when the notion of potential does not make sense).

(i) There is a homomorphism of rings
$$\chi: Mot(X)\to Constr(X,\Z)\,\,,$$ where
$ Constr(X,\Z)$ is the ring of integer-valued constructible functions on $X$ endowed with the pointwise multiplication. More precisely, the element $[\pi:Y\to X]$ is mapped into $\chi(\pi)$, where
$\chi(\pi)(x)=\chi(\pi^{-1}(x))$, which is the Euler characteristic of the fiber $\pi^{-1}(x)$.

(ii) Let now $X$ be a scheme of finite type over a field $\k$, and $l\ne char\,\k$ be a prime number. There is a homomorphism of rings
$$Mot(X)\to K_0(D^b_{\op{constr}}(X,\Q_l))\,\,,$$
 where $D^b_{\op{constr}}(X,\Q_l)$ is the bounded derived category of \'etale $l$-adic  sheaves on  $X$ with constructible cohomology. It is defined by the formula
$$[\pi:S\to X]\mapsto \pi_{!}(\Q_l)\,\,,$$
which is the direct image in the derived sense of the constant sheaf $\Q_l$. Notice that $D^b_{\op{constr}}(X,\Q_l)$ is a tensor category, hence the Grothendieck group $K_0$ is naturally a ring. In this context the homomorphisms $f_{!}$ and $f^{\ast}$ discussed above correspond to the functors $f_{!}$ (direct image with compact support) and $f^{\ast}$ (pullback), which we will denote by the same symbols. We will also use the notation
$\int_X\phi:=f_{!}\,(\phi)$ for the canonical map $f:X\to Spec(\k)$.

(iii) In the special case $X=Spec(\k)$ the above homomorphism becomes a map
$$[S]\mapsto \sum_i(-1)^i[H^i_{c}(S\times_{Spec (\k)} Spec (\overline{\k}),\Q_l)]\in
K_0(\op{Gal}(\overline{\k}/\k)-mod_{\Q_l})\,\,,$$
where $\op{Gal}(\overline{\k}/\k)-mod_{\Q_l}$ is the tensor category of finite-dimensional continuous
 $l$-adic representations of the Galois group $\op{Gal}(\overline{\k}/\k)$, and we take the \'etale cohomology of $S$ with compact support.

(iv) If $\k={\bf F}_q$ is a finite field then for any $n\geqslant 1$ we have a homomorphism
$$Mot(X)\to \Z^{X({\bf F}_{q^n})}$$
 given by $$[\pi:Y\to X]\mapsto (x\mapsto \#\{y\in X({\bf F}_{q^n})\,|\,\pi(y)=x\})\,\,
.$$
Here the operations $f^!,f_\ast$ correspond to  pullbacks and pushforwards of functions
on finite sets.

(v) If $\k\subset \C$ then the category of $l$-adic constructible sheaves on a scheme of finite type $X$ can be replaced in the above considerations by  Saito's category of mixed Hodge modules.

(vi) In the case $X=Spec(\k)$ one has two additional homomorphisms:

a) The Serre polynomial $$Mot(Spec(\k))\to \Z[q^{1/2}]$$ defined by
$$[Y]\mapsto \sum_i(-1)^i\sum_{w\in {\Z_{\geqslant 0}}}\dim H^{i,w}_c(Y)q^{w/2}\,\,,$$
where
$H^{i,w}_c(Y)$ is the weight $w$ component in the $i$-th Weil cohomology group with compact support.

b) If $char\,\k=0$ then we also have the Hodge polynomial $$Mot(Spec(\k))\to \Z[z_1,z_2]$$ given by
$$[Y]\mapsto \sum_{i\geqslant 0}(-1)^i\sum_{p,\,q\geqslant 0}\dim \op{Gr}^p_F(\op{Gr}^W_{p+q} H^i_{DR,c}(Y))z_1^{p}z_2^{q}\,\,,$$
where $\op{Gr}^W_{\bullet} $ and $\op{Gr}_F^{\bullet}$ denote the  graded components with respect to the weight and Hodge filtrations, and $H^i_{DR,c}$ denotes the de Rham cohomology with compact support.

Clearly the Hodge polynomial determines the Serre polynomial via the homomorphism $\Z[z_1,z_2]\to \Z[q^{1/2}]$ such that $z_i\mapsto q^{1/2}, i=1,2$.

Let $X$ be a constructible set over a field $\k$ and $G$ be an affine algebraic group acting on $X$, in the sense that $G(\kk)$ acts on $X(\kk)$ and there
exists a $G$-variety $S$ over $\k$ with a constructible equivariant identification $X(\kk)\simeq S(\kk)$.

We define the group $Mot^G(X)$ of $G$-equivariant motivic functions as the abelian group generated by all $G$-equivariant constructible maps $[Y\to X]$
 modulo the relations
\begin{itemize}
\item $[(Y_1\sqcup Y_2)\to X]=[Y_1\to X]+[Y_2\to X]$,
\item $[Y_2\to X]=[(Y_1\times {\bf A}_\k^d)\to X]$ if $Y_2\to Y_1$ is a $G$-equivariant constructible vector
bundle of rank $d$.
\end{itemize}
This group is a commutative ring via the fiber product, and a morphism of constructible sets with group actions induces a pullback homomorphism
of corresponding rings. There is no natural operation of a pushforward for equivariant motivic functions.

Let $X$ be a constructible set acted by an affine algebraic group $G$. It defines an object $(X,G)$ in the $2$-category of constructible stacks.
The abelian group of stack motivic function $Mot_{st}((X,G))$ is generated by the group of isomorphism classes of 1-morphisms of stacks $[(Y,H)\to (X,G)]$ with the fixed target $(X,G)$, subject to the relations
\begin{itemize}
\item
$[((Y_1,G_1)\sqcup (Y_2,G_2))\to (X,G)]=[(Y_1,G_1)\to (X,G)]+[(Y_2,G_2)\to (X,G)]$
\item $[(Y_2,G_1)\to (X,G)]=[(Y_1\times {\bf A}_\k^d,G_1)\to (X,G)]$ if $Y_2\to Y_1$ is a $G_1$-equivariant constructible vector
bundle of rank $d$.
\end{itemize}
The ring $Mot^G(X)$ maps to $Mot_{st}((X,G))$.

Finally, for a constructible stack ${\cal S}=(X,G)$ we define its class
in the ring $ K_0(Var_{\k})[[{\mathbb L}]^{-1}, ([GL(n)]^{-1})_{n\geqslant 1}]$ as
$$[{\cal S}]={[(X\times GL(n))/G]\over {[GL(n)]}},$$
where we have chosen an embedding $G\to GL(n)$ for some $n\geqslant 1$, and $(X\times GL(n))/G$ is the ordinary quotient by the diagonal free action (thus in the RHS we have the quotient of motives of ordinary varieties). The result does not depend on the choice of embedding. Then we define the integral
$\int_{\cal S}: Mot_{st}({\cal S})\to K_0(Var_{\k})[[{\mathbb L}]^{-1}, ([GL(n)]^{-1})_{n\geqslant 1}]$ as
$\int_{\cal S}[{\cal S}^{\prime}\to {\cal S}]=[{\cal S}^{\prime}]$.

Furthermore, in Section 4.5 of [KoSo1] we associated with the ring of motivic function $Mot(X)$ (and with all equivariant versions of it) another ring $\overline{Mot(X)}$ of certain equivalence classes of motivic functions.  The new ring has e.g. the property that two quadrics which have the same rank and determinant define in it the same element (this is probably not true for the usual motivic rings, see also the section on orientation data below). In Section 4.6 of [KoSo1] we explained that $\overline{Mot(X)}$ and its equivariant cousin can be realized in terms of the rings of functions with numerical (rather than motivic) values.  This property allows us to prove identities in motivic rings by looking at the reduction of the corresponding schemes $mod\,p$ where $p$ is a sufficiently large prime number.

\subsection{The motivic Milnor fiber}

We collect here some facts from [KoSo1], Section 4.3.
Let $M$ be a complex manifold, $x_0\in M$. Recall, that for a germ $f$ of an analytic function at $x_0$ such that $f(x_0)=0$ one can define its {\it Milnor fiber } $MF_{x_0}(f)$, which is a locally trivial $C^{\infty}$-bundle over $S^1$  of manifolds with the boundary (defined only up to a diffeomorphism):
$$ \{z\in M|\,\,dist(z,x_0)\le\varepsilon_1, |f(z)|=\varepsilon_2\}\to S^1=\R/2\pi\Z\,\,,$$
where $z\mapsto \op{Arg}\,f(z)$.
Here $dist$ is any smooth metric on $M$ near $x_0$, and there exists a constant $C=C(f,dist)$ and a positive integer $N=N(f)$ such that for all choices $0<\varepsilon_1\le C$ and $0<\varepsilon_2<\varepsilon_1^N$ the
$C^{\infty}$ type of the bundle is the same for all $\varepsilon_1,\varepsilon_2,dist$.

In particular, taking the cohomology of the fibers of $MF_{x_0}(f)$ we obtain a well-defined local system on $S^1$.

There are several algebro-geometric versions of this construction (theories of nearby cycles). They produce analogs of the local system on $S^1$, for example $l$-adic representations of the group $\op{Gal}(\k((t))^{sep}/\k((t)))$ where $l\ne char\,\k$. In [KoSo1] we used two motivic versions of this notion: one developed by Denef and Loeser and another one (in the framework of non-archimedean analytic geometry) developed by Nicaise and Sebag. In both cases the authors assume that $char\,\k=0$.
Let $X$ be a scheme over $\k$. We introduce the group $\mu=\varprojlim_{n}\mu_n$ and assume that $\mu$ acts trivially on $X$ (here
$\mu_n$ is the group of $n$-th roots of $1$ in $\k$). Then we have the group $Mot^{\mu}(X)$.
We will  assume that the $\mu$-action  is ``good" in the sense that $\mu$ acts via a  finite quotient $\mu_n$ and every orbit is contained in an affine open subscheme.

Let $M$ be a smooth formal scheme over $\k$ with closed point $x_0$ and $f$ be a formal function on $M$
vanishing at $x_0$ (e.g. $M$ could be the formal completion at $0$ of a fiber of vector bundle $V\to X$ in the above notation). We assume that $f$ is not identically equal to zero near $x_0$.

Let us choose a simple normal crossing resolution of singularities $\pi:M'\to M$ of the hypersurface in $M$ given by
the equation $f=0$ with exceptional divisors $D_j, j\in J$. This follows from Hironaka's theorem about resolution of singularities. In fact we need a {\it canonical} resolution of singularities, see e.g. [BiMi], [Te1,2] which better suits the equivariant framework we are dealing with.  Alternatively, as we explained in [KoSo1], one can use the definition of the motivic Milnor fiber based on Berkovich theory of non-archimedean analytic spaces (see [NiSe]). Both approaches work in case when $f$ is a formal series (e.g. in the case of the approach via canonical resolution of singularities it follows from  [Te3]).
The explicit formula for the motivic Milnor fiber  in terms of a resolution of singularities looks as follows

$$MF_{x_0}(f)=\sum_{I\subset J, I\ne \emptyset}(1-{\mathbb L})^{\#I-1}[\widetilde{{D_I}^0} \cap \pi^{-1}(x_0)]\in Mot^\mu
(Spec(\k))\,\,,$$
where $D_I=\cap_{j\in I}D_j$, $D_I^0$ is the complement in $D_I$ to the union of all other exceptional divisors,
and $\widetilde{{D_I}^0}\to D_I^0$ is a certain Galois cover with  Galois group  $\mu_{m_I}$, where $m_I$ is the g.c.d.
of the multiplicities of all divisors $D_i,i\in I$. Informally speaking, the fiber of the cover
 $\widetilde{{D_I}^0}\to D_I^0$ is the set of connected components of a non-zero level set of the function $f\circ \pi$ near a point of  $D_I^0$.

Furthermore, the group $Mot^{\mu}(X)$ carries a non-trivial associative commutative product introduced by Looijenga (which we call ``exotic product" in Section 4.3 of [KoSo1]). We define
$${\cal M}^{\mu}(X):=(Mot^{\mu}(X), \mathrm{\, exotic\,\,\,product})\,\,. $$

Let $V\to X,\,V'\to Y$ be two constructible vector bundles endowed with constructible families $f,g$ of formal power series.
We denote by $f\oplus g$ the sum of pullbacks of $f$ and $g$ to the constructible vector
bundle
$$pr_X^*V\oplus pr_Y^*V'\to X\times Y\,\,.$$
The main result on motivic Milnor fibers that we need is the following motivic Thom-Sebastiani theorem proved by Denef and Loeser (see [DeLo]).

\begin{thm} One has
$$(1-MF(f\oplus g))=pr_X^*(1-MF(f))\cdot
pr_Y^*(1-MF(g))\in {\cal M}^\mu (X\times Y)\,\,.$$
\end{thm}

The version of the motivic Milnor fiber defined by means of non-archimedean analytic geometry in [NiSe] agrees with the formula of Denef and Loeser. On the other hand, in the non-archimedean approach one can mimick the classical construction which we recalled at the beginning of this section.
Also, the above definitions can be extended  to the equivariant setting as well as to the case when $f$ depends on parameters. For our purposes the function $f$ will be the potential $W_E$ considered in a neighborhood of an object.

\subsection{Orientation data}

The definition of motivic DT-invariants depends on the minimal model of the potential. This means that we make a non-commutative change of formal variables so that the Taylor series of the potential $W_E$ starts with the terms of degree at least $3$. This property is not invariant under the change of an object: quadratic terms can appear if we decompose the same series at a nearby object. It turns out that the consistent choice of the quadratic part depends on an additional structure, which we call {\it orientation data} in [KoSo1], Section 5. We recall it here.

Let $\CC$ be an ind-constructible $\k$-linear $3$-dimensional Calabi-Yau category.
Then we have a natural ind-constructible super line bundle ${\cal D}$ over $Ob(\CC)$ with the fiber over $E$ given by
${\cal D}_E=\op{sdet}(\op{Ext}^{\bullet}(E,E))$. It follows that on the ind-constructible stack of exact triangles $E_1\to E_2\to E_3$ we have an isomorphism of the pull-backed line bundles which fiberwise reads as
$${\cal D}_{E_2}\otimes {\cal D}_{E_1}^{-1}\otimes {\cal D}_{E_3}^{-1}\simeq (\op{sdet}(\op{Ext}^{\bullet}(E_1,E_3)))^{\otimes 2}\,\,.$$

Let $Ob(\CC)=\sqcup_{i\in I}Y_i$ be a decomposition into a disjoint countable union of constructible sets,
each acted by an affine algebraic group $GL(N_i)$.

\begin{defn} Orientation data on $\CC$ consists of a choice of an ind-constructible super line bundle $\sqrt{{\cal D}}$ on $Ob(\CC)$ such that its restriction to each $Y_i, i\in I$ is $GL(N_i)$-equivariant, endowed on each $Y_i$ with $GL(N_i)$-equivariant isomorphisms $(\sqrt{{\cal D}})^{\otimes 2}\simeq {\cal D}$ and such that for the natural pull-backs to the ind-constructible stack of exact triangles $E_1\to E_2\to E_3$ we are given equivariant isomorphisms:
$$\sqrt{{\cal D}}_{E_2}\otimes (\sqrt{{\cal D}}_{E_1})^{-1}\otimes (\sqrt{{\cal D}}_{E_3})^{-1}\simeq \op{sdet}(\op{Ext}^{\bullet}(E_1,E_3))$$
such that the induced equivariant isomorphism
$${\cal D}_{E_2}\otimes {\cal D}_{E_1}^{-1}\otimes {\cal D}_{E_3}^{-1}\simeq (\op{sdet}(\op{Ext}^{\bullet}(E_1,E_3)))^{\otimes 2}$$
coincides  with the one which we have a priori.
\end{defn}

Orientation data are in a sense similar to a choice of spin structure on Lagrangian submanifolds in the definition of Fukaya category (there are several other situations where one has to make such a choice in order to define objects of a category). As we mentioned above, our story  is related to the fact that the minimal model of the potential $W_E$ at an object $E$ gives rise to a non-minimal potential $W_F$ for a close object $F$.
To say it differently, let $V$ be a $\k$-vector space endowed with a non-degenerate quadratic form $Q$. We define an element
$${I}(Q)=(1-MF_0(Q)){\mathbb L}^{-{1\over{2}}\dim V}\in {\cal M}^{\mu}(Spec(\k))[{\mathbb L}^{\pm 1/2}]\,\,,$$
 where
${\mathbb L}^{1/2}$ is a formal symbol which satisfies the relation $({\mathbb L}^{1/2})^2={\mathbb L}$, and $Q$ is interpreted as a function on $V$.
Then the motivic Thom-Sebastiani theorem implies
that $${I}(Q_1\oplus Q_2)={I}(Q_1){I}(Q_2)\,\,.$$
 Also we have ${I}(Q)=1$, if $Q$ is a split form: $Q=\sum_{1\le i\le n}x_iy_i$
for $V=\k^{2n}$. This means that we can ``twist" the element $(1-MF(W))$ by $I(Q)$.
In all motivic realizations (e.g. by taking the Serre polynomial) the element $I(Q)$ depends only on the pair
$(\op{rk}Q, \op{det}(Q)\op{mod}(\k^{\times})^2)$. It is not known (and probably not true) if the motive of a quadric
is uniquely determined by the above pair of numbers. This problem readily translates to the case of constructible families. As we will see below, the homomorphism from the motivic Hall algebra to the motivic quantum torus depends on the choice of quadratic part $Q_E$ of the potential $W_E$ (or rather on the corresponding element $I(Q)$). In order to make a consistent choice of the latter we need to pick orientation data. The general problem of its existence is open (the answer is positive for $3CY$ categories corresponding to quivers with potential).

\section{Motivic Hall algebra}

Let $Ob(\CC)=\sqcup_{i\in I}X_i$ be a constructible decomposition of the space of objects with the group $GL(N_i)$ acting on $X_i$.
Let us consider the $Mot(Spec(\k))$-module $\oplus_i Mot_{st}(X_i,GL(N_i))$ and extend it by adding free generators which are negative powers ${\mathbb L}^n, n<0$ of the motive of the affine line  ${\mathbb L}$. We denote the resulting $Mot(Spec(\k))$-module by $H(\CC)$.  We understand elements of
$H(\CC)$ as {\it measures} (and not as functions), because in the definition of the product we will use the pushforward maps.
In the Section 6.1 of [KoSo1] we introduced a structure of an asociative algebra on $H(\CC)$ and called it {\it motivic Hall algebra}. The product is defined on constructible families by the formula
$$[Y_1\to Ob(\CC)]\cdot [Y_2\to Ob(\CC)]=\sum_{n\in \Z}[W_n\to Ob(\CC)]{\mathbb L}^{-n}\,\,,$$
where
$$W_n  =\big\{(y_1,y_2,\alpha)\,|\,y_i\in Y_i, \alpha\in \op{Ext}^1(\pi_2(y_2),\pi_1(y_1))\,,\,
 (\pi_2(y_2),\pi_1(y_1))_{\le 0}=n\big\}\,\,, $$
with the notation
 $$(E,F)_{\le N}:=\sum_{i\le N}(-1)^i\dim \op{Ext}^i(E,F)\,\,.$$
The map $W_n\to Ob(\CC)$ is given by the formula
$$(y_1,y_2,\alpha)\mapsto Cone(\alpha: \pi_2(y_2)[-1]\to \pi_1(y_1))\,\,.$$

In other words, the structure constants are given by motives of exact triangles.
For a constructible stability structure on $\CC$ with an ind-constructible class map $\op{cl}:K_0(\CC)\to \Gamma$, a central charge $Z:\Gamma\to \C$, a strict sector $V\subset \R^2$ and a branch $\op{Log}$ of the logarithm function on $V$ we have the category $\CC_V:=\CC_{V,\op{Log}}$ with objects which are either the zero object or an extension of semistable objects $E$ such that $Z(cl(E))\in V$. Hence we have
the associative algebra given by the {\it completion}
$$\widehat{H}(\CC_V):=\prod_{\gamma\in (\Gamma\cap C(V,Z,Q))\cup\{0\}}H(\CC_V\cap \op{cl}^{-1}(\gamma))\,\,,$$
where $C(V,Z,Q)$ is the convex cone spanned by elements $\gamma\ne 0$ such that $Z(\gamma)\in V$ and $Q(\gamma)\ge 0$ (recall here the quadratic form $Q$ introduced at the end of Section 2.1).

For every such algebra  we define an invertible element $A_V^{\op{Hall}}\in \widehat{H}(\CC_V)$ such that
$$A_V^{\op{Hall}}:=1+\dots=\sum_{i\in I}  \mathbf{1}_{(Ob(\CC_V)\cap Y_i,GL(N_i))}\,\,,$$
where $\mathbf{1}_{\cal S}$ is the identity function but interpreted as a counting measure. Then we proved the following easy result.

\begin{prp} The elements $A_V^{\op{Hall}}$ satisfy the Factorization Property:
$$A_{V}^{\op{Hall}}=A_{V_1}^{\op{Hall}}\cdot A_{V_2}^{\op{Hall}}$$
for a strict sector $V=V_1\sqcup V_2$ (decomposition in the clockwise order).
\end{prp}

\section{Quantum torus and motivic DT-invariants}

\subsection{Motivic weights}

Let $\CC$ be a $3$-dimensional ind-constructible   Calabi-Yau category over a field $\k$ of characteristic zero.
We may assume that it is {\it minimal on the diagonal}, which means that we replaced every $\A$-algebra $End(E,E)$ by its minimal model (i.e. $m_1=0$).
Then for any $E\in Ob(\CC)(\kk)$ we have  the potential $W_E^{min}$ which is a formal power series in $\alpha\in \op{Ext}^1(E,E)$ which starts with cubic terms. It is just the minimal model of $W_E$.
We denote by
$$MF(E):=MF_0(W_E^{min})$$ the motivic Milnor fiber of $W_E^{min}$ at $0\in \op{Ext}^1(E,E)$. Then the assignment $E\mapsto MF(E)$ can be interpreted as the value of some function $MF\in {\cal M}^{\mu}(Ob(\CC))$, where $\mu$, as before, is the group of all roots of $1$.

Let us choose  orientation data $\sqrt{\cal D}$ for $\CC$. This allows us to speak about a super line bundle $\sqrt{\cal D}\otimes {\cal D}_{\le 1}^{-1}$ with trivialized tensor square. Here
${\cal D}_{\le 1}$ is a super line bundle over the space of objects with the  fiber at $E\in Ob(\CC)$ given by
$${\cal D}_{\le 1,E}:=\op{sdet}(\tau_{\le 1}(\op{Ext}^{\bullet}(E,E)))\,\,,$$ where $\tau_{\le i}, i\in \Z$ denotes the standard truncation functor.

The data of a super line bundle $\sqrt{\cal D}\otimes {\cal D}_{\le 1}^{-1}$ is basically the same the data consisting of a pair $(V,Q)$ where $V$ is a super vector space endowed with a quadratic form $Q$. Thus we can define
${I}(Q)=(1-MF_0(Q)){\mathbb L}^{-{1\over{2}}\op{rk} Q}$ as an element of the appropriate motivic ring.

\begin{defn} The motivic weight $w\in {\cal M}^{\mu}(Ob(\CC))$ is the function defined on objects by the formula
$$w(E)={\mathbb L}^{{1\over{2}}\sum_{i\le 1}(-1)^i \dim \op{Ext}^i(E,E)}(1-MF(E))(1-MF_0(Q_E)){\mathbb L}^{{-{1\over{2}}\op{rk}Q_E}}\,\,.$$

\end{defn}

\subsection{The motivic quantum torus}

Let $\Gamma$ be a free abelian group endowed with a skew-symmetric, integer-valued bilinear form $\langle\bullet,\bullet\rangle$.
For any commutative unital ring $C$ which contains an invertible symbol ${\mathbb L}^{1/2}$ we introduce a $C$-linear associative algebra
$${\cal R}_{\Gamma,C}:=\oplus_{\gamma\in \Gamma}C\cdot \hat{e}_{\gamma}$$
 where the generators $\hat{e}_{\gamma}, \gamma\in \Gamma$ satisfy the relations
$$\hat{e}_{\gamma_1}\hat{e}_{\gamma_2}={\mathbb L}^{{1\over{2}}\langle \gamma_1,\gamma_2\rangle}\hat{e}_{\gamma_1+\gamma_2}, \,\,\,\hat{e}_0=1\,\,.$$
 We will call it the  {\it quantum torus} associated with $\Gamma$ and $C$.

For any strict sector $V\subset \R^2$ we define
$${\cal R}_{V,C}:=\prod_{\gamma \in \Gamma\cap \,C_0(V,Z,Q)}C\cdot \hat{e}_{\gamma}$$
and call it the quantum torus associated with $V$. Here
$$C_0(V,Z,Q):=C(V,Z,Q)\cup \{0\}$$
and $C(V,Z,Q)$  is the convex cone generated by the set of vectors $x\in \Gamma_{\R}:=\Gamma\otimes \R\setminus \{0\}$
such that $Z(x)\in V$ and $Q(x)\ge 0$.
The algebra  ${\cal R}_{V,C}$ is the natural completion of the subalgebra  ${\cal R}_{V,C}\cap {\cal R}_{\Gamma,C}\subset{\cal R}_{\Gamma,C}$. Notice that a choice of the sector $V$ allows us to define  the completion of ${\cal R}_{\Gamma,C}$ which contains the motivic DT-invariant (see below). Being an infinite series the motivic DT-invariant does not belong to ${\cal R}_{\Gamma,C}$.

Consider next the ring
$${D}^{\mu}={\cal M}^{\mu}(Spec(\k))[{\mathbb L}^{-1}, {\mathbb L}^{1/2},  ([GL(n)]^{-1})_{n\geqslant 1}]\,\,, $$
where ${\mathbb L}=[{\bf A}^1_{\k}]$ is the motive of the affine line. The element
${\mathbb L}^{1/2}$ is a formal symbol satisfying the equation $({\mathbb L}^{1/2})^2={\mathbb L}$. Instead of inverting the motives $$[GL(n)]=({\mathbb L}^n-1)({\mathbb L}^n-{\mathbb L})\dots({\mathbb L}^n-{\mathbb L}^{n-1})$$ of all general linear groups we can invert motives of all projective spaces
$$[{\bf P}^n]={{\mathbb L}^{n+1}-1\over {\mathbb L}-1}\,\,.$$
We also will consider the ring $\overline{D^{\mu}}$ of equivalence classes of functions in $D^{\mu}$ mentioned
in the section 2.2.
The ring $\overline{D^{\mu}}$ will play the role of the universal coefficient ring where motivic Donaldson-Thomas invariants take value. In particular, we choose it as the coefficient ring $C$ in the above definition of the quantum torus.
We denote by ${\cal R}_{\Gamma}:={\cal R}_{\Gamma, \overline{D^{\mu}}}$
the corresponding quantum torus and call it the {\it motivic quantum torus} associated with $\Gamma$. Similarly, we have  motivic quantum tori ${\cal R}_V$ associated with strict sectors $V\subset \R^2$.

\subsection{From motivic Hall algebra to motivic quantum torus}

Assume  that $\CC$ is an ind-constructible  $3d$ Calabi-Yau category over a field $\k$ of characteristic zero, endowed with a constructible stability condition and orientation data $\sqrt{\cal D}$. The Hall algebra of $\CC$ is graded by the corresponding lattice $\Gamma$: $H(\CC)=\oplus_{\gamma\in \Gamma}H(\CC)_{\gamma}$.  The following theorem was formulated in Section 6.3 of [KoSo1].

\begin{thm} The map $\Phi: H(\CC)\to {\cal R}_{\Gamma}$ given by the formula
$$\Phi(\nu)=(\nu,w)\hat{e}_{\gamma}, \,\,\,\nu\in H(\CC)_{\gamma}$$
is a homomorphism of $\Gamma$-graded $\Q$-algebras. Here $w$ is the motivic weight and $(\bullet,\bullet)$ is the pairing between motivic measures and functions.

\end{thm}

In other words, the homomorphism $H(\CC)\to {\cal R}_{\Gamma}$ can be written as
$$ [\pi: Y\to Ob(\CC)]\mapsto $$
$$\mapsto \int_{Y}(1-MF(\pi(y)))\,(1-MF_0(Q_{\pi(y)}))\,{\mathbb L}^{-{1\over{2}}\op{rk}Q_{\pi(y)}}\,{\mathbb L}^{{1\over{2}}(\pi(y),\pi(y))_{\le 1}}\,\hat{e}_{\op{cl}(\pi(y))}\,\,,$$
where $\int_Y$ is understood as the direct image functor.

The natural extension of the above homomorphism to the completion of $\widehat{H}(\CC_V)$  maps the element $A_V^{\op{Hall}}$ to the
element $A_V^{mot}$.

The proof of the above theorem given in [KoSo1] is based on the motivic Thom-Sebastiani theorem and a certain identity of motivic integrals formulated in Section 4.4 of loc. cit. We do not reproduce the identity here. We formulated it as the Conjecture 4 and gave a proof of it in the $l$-adic case (i.e. when instead of the theory of motivic functions and its equivariant version one uses the corresponding well-known theory of constructible $l$-adic functions). Therefore, technically speaking, the above theorem was proved in [KoSo1] only for the $l$-adic quantum torus and the $l$-adic Hall algebra rather then for their motivic cousins. The $l$-adic versions however suffice if we are interested in numerical DT-invariants. Indeed, the results of Section 4.6 of [KoSo1] imply that one can use the $l$-adic realization of motivic DT-invariants in order to obtain the numerical DT-invariants (modulo the existence of the quasi-classical limit, as we explain below).
In order to simplify the notation and exposition we are going to assume below that  Theorem 2 is true in its full generality (i.e. that the Conjecture 4 from [KoSo1] holds not only $l$-adically, but for motives). Then we can define the collection $\Phi(A_V^{Hall}):=A_V^{mot}$ parametrized by strict sectors $V\subset \R^2$.

\begin{defn} Let $\CC$ be an ind-constructible $3$-dimensional Calabi-Yau category endowed with a stability condition $\sigma\in Stab(\CC,\op{cl})$  and an orientation data. We call the collection of elements $(A_V^{mot}\in {\cal R}_V)$ of the completed motivic quantum tori $({\cal R}_V)$ (for all strict sectors $V\subset \R^2$) the  motivic Donaldson-Thomas invariant of $\CC$.

\end{defn}

\subsection{Quasi-classical limit and integrality}

Let us consider the following unital $\Q$-subalgebra of $\Q(q^{1/2})$:
$$D_q:=\Z[q^{1/2},q^{-1/2},\left((q^n-1)^{-1}\right)_{n\geqslant 1}]\,.$$
 There is a homomorphism of rings $\overline{D^{\mu}}\to D_q$ given by the twisted Serre polynomial. Namely, it maps ${\mathbb L}^{1/2}\mapsto q^{1/2}$, and on
${\cal M}^{\mu}$ it is the composition of the Serre polynomial with the involution $q^{1/2}\mapsto -q^{1/2}$.
Therefore, we have a homomorphism of algebras ${\cal R}_{\Gamma}\to {\cal R}_{\Gamma,q}$, where
${\cal R}_{\Gamma,q}$ is the $D_q$-algebra
generated by $\hat{e}_{\gamma}, \gamma\in \Gamma$, subject to the relations
$$\hat{e}_{\gamma}\hat{e}_{\mu}=q^{{1\over{2}}{\langle\gamma,\mu\rangle}}\hat{e}_{\gamma+\mu}\,,\,\,\,\hat{e}_0=1\,.$$

Here the integer skew-symmetric form $\langle\bullet,\bullet\rangle$ on $\Gamma$ is a part of the data defining a stability condition on $\CC$ (see Introduction in [KoSo1] for the details and explanations on why our axioms are slightly different from those of Bridgeland).
Similarly to the motivic case, we have an algebra ${\cal R}_{V,q}$ associated with any strict sector $V$.

The elements $A_{V,q}\in {\cal R}_{V,q}$ corresponding to $A_V^{mot}$ are series in $\hat{e}_{\gamma},\gamma\in \Gamma$ with coefficients which are rational functions in $q^{1/2}$. They can have poles as $q^n=1$ for some $n\geqslant 1$. Hence it is not clear how to take the quasi-classical limit as $q^{1/2}\to -1$ (this corresponds to the taking of Euler characteristic of the corresponding motives).

Let us assume that the skew-symmetric form on $\Gamma$ is non-degenerate (otherwise we can replace $\Gamma$ by the symplectic lattice $\Gamma\oplus \Gamma^{\vee}$). The element $A_{V,q}$ defines an automorphism of an appropriate completion of
${\cal R}_{\Gamma,q}$. More precisely, it acts by the conjugation $x\mapsto A_{V,q}xA_{V,q}^{-1}$ on the subring
$$\prod_{\gamma\in C_0(V)\cap \Gamma}D_q\hat{e}_{\gamma}\,$$ where $C_0(V)=C_0(V,Z,Q)$ is the union of $0$ with the
convex hull $C(V,Z,Q)$ of the set $Z^{-1}(V)\cap \{Q\geqslant 0\}$.

 The ``integer" quantum torus
$$\bigoplus_{\gamma\in C_0(V)\cap \Gamma} \Z[q^{\pm 1/2}]\hat{e}_{\gamma}\subset {\cal R}_{\Gamma,q}$$ has the quasi-classical limit\footnote{There is another quasi-classical limit $q^{1/2}\to +1$ which we do not consider here.} which is the Poisson algebra with basis ${e}_{\gamma}, \gamma\in C_0(V)\cap \Gamma$ with  product and Poisson bracket
given by
$${e}_{\gamma}{e}_{\mu}=(- 1)^{\langle \gamma,\mu\rangle}e_{\gamma+\mu},\,\,\,\{e_{\gamma},e_{\mu}\}=(- 1)^{\langle \gamma,\mu\rangle}\langle \gamma,\mu\rangle e_{\gamma+\mu}\,\,.$$
The Poisson bracket is the limit of a normalized bracket:
$$[\hat{e}_\gamma,\hat{e}_\mu]=\left(q^{1/2\langle \gamma, \mu\rangle}-q^{-1/2\langle \gamma, \mu\rangle}\right)\hat{e}_{\gamma+\mu}\,\,,$$
$$\lim_{q^{1/2}\to -1}
(q-1)^{-1}\cdot\left(q^{1/2\langle \gamma, \mu\rangle}-q^{-1/2\langle \gamma, \mu\rangle}\right)=(- 1)^{\langle \gamma,\mu\rangle}\langle \gamma,\mu\rangle\,.$$
 One can write informally
 $$e_\gamma=\lim_{q^{1/2}\to -1}\frac{\hat{e}_\gamma}{q-1}\,\,.$$

 The following ``absense of poles conjecture" was formulated in [KoSo1], Section 7.1.

\begin{conj} For any $3d$ Calabi-Yau category with a stability condition and any strict sector $V$ the automorphism $x\mapsto A_{V,q}xA_{V,q}^{-1}$ preserves the subring $$\prod_{\gamma\in C_0(V)\cap \Gamma}D_q^+\hat{e}_{\gamma}\,\,,$$ where
$D_q^+ :=\Z[q^{\pm 1/2}]$.

\end{conj}

Furthermore, in Section 7.3 of [KoSo1] we formulated a stronger conjecture at the level of motivic Hall algebras. It seems that the above ``absense of poles conjecture" holds for a larger class of categories than  the class of $3CY$ categories.

Assuming the above conjecture we  denote by $A_V$ the ``quasi-classical limit" of the automorphism $x\mapsto A_{V,q}xA_{V,q}^{-1}$ as $q^{1/2}\to -1$.
This is a symplectomorphism of the symplectic torus associated with the symplectic lattice $\Gamma$ (or, better,
with the symplectic lattice $\Gamma\oplus \Gamma^{\vee}$). This torus can be thought of as the quasi-classical limit of the motivic quantum torus. For a generic stability condition this symplectomorphism can be written as

$$A_V=\prod_{Z(\gamma)\in V}^{\longrightarrow}T_{\gamma}^{\,\Omega(\gamma)},$$
where $$T_{\gamma}(e_{\mu})=(1-e_{\gamma})^{\langle \gamma,\mu\rangle}e_{\mu}$$ and $\Omega(\gamma)\in \Q$.
The collection $(A_V)$ satisfies the Factorization Property.
In the same Section 7.1 of [KoSo1] we formulated the following conjecture.
\begin{conj} For a generic central charge $Z$ all numbers $\Omega(\gamma),\gamma \in \Gamma\setminus \{0\}$ are integers.

\end{conj}

We discussed non-trivial arguments in favor of this ``integrality conjecture" in Section 7.5 of [KoSo1]. In the case of the category generated by one Schur object $E$ (i.e. $Ext^{<0}(E,E)=0$ and $Hom(E,E)$ is one-dimensional) we reformulated it in terms of the generating function of Euler characteristics of certain moduli spaces of representations of a quiver with potential. Although in general the conjecture is still open, there are some partial results confirming it. In particular, Reineke proved in [Re1] that the integrality property is preserved under the wall-crossing formula (see below). In [KoSo3] we obtained a stronger than [Re1]  integrality result at the level of motivic DT-invariants, including the integrality in case of quivers with potential.

\begin{rmk}
The collection $(\Omega(\gamma))_{\gamma \in \Gamma}$
seems to be the correct mathematical definition of the counting of BPS states in  String Theory.
\end{rmk}

The quasi-classical limit of motivic DT-invariants can be compared with the ``microlocal" version of DT-invariants introduced by Kai Behrend in [Be1].
Recall that he defined a $\Z$-valued invariant of a critical point $x$ of a function $f$ on $X$ which is equal to
$$(-1)^{\dim X}(1-\chi(MF_x(f)))\,\,,$$ where $\chi$ denotes the Euler characteristic. By Thom-Sebastiani theorem this number does not change if we add to $f$ a function with a quadratic singularity at $x$ (stable equivalence).

Let $M$ be a scheme with symmetric perfect obstruction theory (see [BeFa]). Then $M$ is locally represented as a scheme of critical points of a function $f$ on a manifold $X$. The above invariant gives rise to a $\Z$-valued constructible function $B$ on $M$. The
global invariant is
$$\int_MB\,d\chi:=\sum_{n\in \Z}n \chi(B^{-1}(n))\,,$$ where $\chi$ denotes the Euler characteristic. Behrend proved that for a {\it proper} $M$ the invariant $\int_M B\, d\chi$ coincides with the degree of the virtual fundamental class $[M]^{virt}\in H_0(M)$ given by $\int_{[M]^{virt}}1$.

Now let us assume that $M\subset \CC^{ss}$ consists of Schur objects. For such an object $E$ let $\op{cl}(E)=\gamma\in \Gamma$ be a fixed primitive class.
Then the contribution of $M$ into the $\gamma$-component of the motivic DT-invariant is equal to
$$\int_M{{\mathbb L}^{{1\over{2}}(1-\dim \op{Ext}^1(E,E))}\over {\mathbb L}-1}(1-MF(E))(1-MF_0(Q_E)){\mathbb L}^{-{1\over{2}}\op{rk}Q_E}\hat{e}_{\gamma}\,\,.$$
In the quasi-classical limit this contribution is easily seen to be equal to $-\Omega(\gamma)$.
Hence Behrend formula implies that the contribution of $M$ to the value
$\Omega(\gamma)$ is equal to $\int_{[M]^{virt}}1$.

\begin{rmk} An important property of any generalization of DT-invariants is their deformation invariance. In Section 7.2 of [KoSo1] we formulated a conjecture which claims the invariance of the collection
$(\Omega(\gamma))_{\gamma \in \Gamma}$ with respect to the ``polarization preserving" deformations of $\CC$, in the case when $\CC$ is homologically smooth (the latter is the categorical analog of proper smooth schemes introduced in [KoSo2]). In the geometric framework (i.e. for the category of coherent sheaves) there is a theory of generalized DT-invariants developed by Joyce and Song (see [JoS1,2]). In this more restrictive setting they can prove deformation invariance property of their invariants.

\end{rmk}

\section{Stability data on graded Lie algebras and wall-crossing formulas}

In Section 2 of [KoSo1] we developed a general approach to ``numerical" DT invariants, which is independent on the Calabi-Yau property and even of  categories. It is based on the notion of stability data on graded Lie algebras.

\subsection{Stability data}

Let us fix a free abelian group $\Gamma$ of finite rank,
and a graded Lie algebra $\g=\oplus_{\gamma\in \Gamma}\g_{\gamma}$ over $\Q$.\footnote{In examples $\g$ is a $R$-linear Lie algebra, where $R$ is a commutative unital $\Q$-algebra.}

\begin{defn} A stability data on $\g$ is a pair
$\sigma=(Z,a)$ such that:

1) $Z: \Gamma\to \R^2\simeq \C$ is a homomorphism of abelian groups called the central charge;

2) $a=(a(\gamma))_{\gamma\in \Gamma\setminus \{0\}}$ is a collection of elements $a(\gamma)\in \g_{\gamma}$,

satisfying the following

\vspace{2mm}

{\bf Support Property}:
\vspace{2mm}

Pick a  norm $\parallel\bullet\parallel$ on the vector space $\Gamma_{\R}=\Gamma\otimes_{\Z}\R$.
Then there exists $C>0$ such that for any $\gamma\in \op{Supp}a $ (i.e. $a(\gamma)\ne 0$) one has
$$\parallel\gamma\parallel\le C|Z(\gamma)|\,\,.$$

\end{defn}

Obviously the Support Property does not depend on the choice of the norm. We will denote the set of all stability
data on $\g$ by $Stab(\g)$. Later we will equip this set with a Hausdorff topology.

The Support Property is equivalent to the following condition (which we will also call the Support Property):

{\it There exists a  quadratic form $Q$ on $\Gamma_{\R}$ such that

1) $Q_{|\op{Ker}\, Z}<0$;

2) $\op{Supp} a\subset \{\gamma\in \Gamma\setminus \{0\}|\,\,Q(\gamma)\geqslant 0\}$,

where we use the same notation $Z$ for the natural extension of $Z$ to $\Gamma_{\R}$.}

Indeed, we may assume that the norm $\parallel\bullet\parallel$ is the Euclidean norm in a chosen basis and
take $Q(\gamma)=-\parallel\gamma\parallel^2+C_1|Z(\gamma)|^2$ for sufficiently large positive
constant $C_1$. Generically $Q$ has signature $(2,n-2)$, where $n=\op{rk}\Gamma$.

For a given quadratic form $Q$ on $\Gamma_{\R}$ we denote by $Stab_Q(\g)\subset Stab(\g)$ the set of stability data satisfying the above conditions 1) and 2). Obviously $Stab(\g)=\cup_QStab_Q(\g)$, where the union is taken over all quadratic forms $Q$.

In Section 2.2 of [KoSo1] we reformulated the stability data in the following equivalent form.

We denote by $\widehat{Stab}_Q(\g)$ the set of pairs $(Z, A)$ such that:

{\it a) $Z: \Gamma\to \R^2$ is an additive map such that
$Q_{|\op{Ker}Z}<0$;

b) $A=(A_V)_{V\in {\cal S}}$ is a collections of elements $A_V\in G_{V,Z,Q}$, where
$G_{V,Z,Q}$ is a pronilpotent
group with the pronilpotent graded Lie algebra
$$\g_{V,Z,Q}=\prod_{\gamma\in\Gamma\cap C(V,Z,Q)}\g_{\gamma}\,\,,$$
where $C(V,Z,Q)$ is the convex cone generated by the set
$$S(V,Z,Q)=\{x\in \Gamma_{\R}\setminus \{0\}|\, Z(x)\in V, Q(x)\geqslant 0\}\,\,.$$}

For a triangle $\Delta$ which is cut out of the sector $V$ by a straight line,  any $\gamma\in Z^{-1}(\Delta)$ can be represented as a sum of  elements
of $\Gamma\cap C(V,Z,Q)$ in finitely many ways.
Furthermore, the triangle $\Delta$ defines the Lie ideal $J_{\Delta}\subset \g_{V,Z,Q}$ consisting
of elements $y=(y_{\gamma})\in \g_{V,Z,Q}$ such that for every component $y_{\gamma}$
the corresponding $\gamma$ does not belong to the convex hull
of $Z^{-1}(\Delta)$. Then the quotient $\g_{\Delta}:=\g_{V,Z,Q}/J_{\Delta}$ is a  nilpotent
Lie algebra, and $\g_{V,Z,Q}=\varprojlim_{\Delta\subset V}g_{\Delta}$.
Let $G_{\Delta}=\exp(\g_{\Delta})$ be the pronilpotent  group with the pronilpotent Lie algebra $\g_{\Delta}$.
Then $G_{V,Z,Q}=\varprojlim_{\Delta}G_{\Delta}$ is a pronilpotent group.
If $V=V_1\sqcup V_2$ (in the clockwise order) then there are natural embeddings
$G_{V_i,Z,Q}\to G_{V,Z,Q}, i=1,2$.
We impose the following axiom on the set of pairs $(Z,A)$:

{\bf Factorization Property:}

{\it The element $A_V$ is given by the product
$A_V=A_{V_1}A_{V_2}$ where the equality is understood in $G_{V,Z,Q}$.}

If we assume the ``quasi-classical limit" conjecture discussed in Section 4.4 then the algebra of functions on the  symplectic torus $(\Gamma\oplus \Gamma^{\vee})\otimes \R/\Gamma\oplus \Gamma^{\vee}$ gives rise to a $\Gamma\oplus \Gamma^{\vee}$-graded Lie algebra. The collection of symplectomorphisms $(A_V)$ obtained  as the limit of the automorphisms $(Ad(A_{V,q}))$ as $q^{1/2}\to -1$ gives rise to  stability data on this graded Lie algebra.

\subsection{Wall-crossing formulas}

In Section 2.3 of [KoSo1] we defined a Hausdorff topology on the space $Stab(\g)$.
Let us recall this definition here.
Let $X$ be a topological space, $x_0\in X$ be a point, and $(Z_x,a_x)\in Stab(\g)$ be a family parametrized by $X$.

\begin{defn} We say that the family is continuous at $x_0$ if the following conditions are satisfied:

a) The map $X\to \op{Hom}(\Gamma,\C), x\mapsto Z_x$ is continuous at $x=x_0$.

b) If a quadratic form $Q_0$ such that
$(Z_{x_0},a_{x_0})\in Stab_{Q_0}(\g)$ then there exists an open
neighborhood $U_0$ of $x_0$ such that
$(Z_{x},a_{x})\in Stab_{Q_0}(\g)$ for all $x\in U_0$.

c) For any closed strict sector $V$ such that $Z(\op{Supp}a_{x_0})\cap \partial V=\emptyset$ the map
$$x\mapsto \log\,A_{V,x,Q_x}\in \g_{V,Z_x,Q_x}\subset \prod_{\gamma\in \Gamma}\g_{\gamma}$$
 is continuous
at $x=x_0$. Here we endow the vector space $\prod_{\gamma\in \Gamma}\g_{\gamma}$ with the product topology of discrete sets, and $A_{V,x,Q_x}$ is the group element associated with $(Z_x,a_x)$, the sector $V$ and a quadratic form $Q_x$ such that
$(Z_x,a_x)\in Stab_{Q_x}(\g)$.

\end{defn}

If $\g$ arises from a $3CY$ category then this definition has a nice categorical interpretation given in Section 3.4 of loc. cit. Having a topology defined in such a way one immediately arrives at the wall-crossing formulas for the stability data such as follows.

Let us fix an element $Z_0\in \op{Hom}(\Gamma,\C)$ and a quadratic form $Q_0$ compatible with $Z_0$ (i.e. negative on its kernel).
We denote by $U_{Q_0,Z_0}$ the connected component containing $Z_0$ in the domain $\{Z\in \op{Hom}(\Gamma,\C)|\,\, (Q_0)_{|\op{Ker} Z}<0\}$.
Let $\gamma_1,\gamma_2\in \Gamma\setminus \{0\}$ be two $\Q$-linearly independent elements such that
$Q_0(\gamma_i)\geqslant 0, Q_0(\gamma_1+\gamma_2)\geqslant 0, i=1,2$.
We introduce the set
$${\cal W}^{Q_0}_{\gamma_1,\gamma_2}=\{Z\in U_{Q_0,Z_0}|\,\,\R_{>0}Z(\gamma_1)=\R_{>0}Z(\gamma_2)\}\,\,.$$
In this way we obtain a countable collection of real hypersurfaces
${\cal W}^{Q_0}_{\gamma_1,\gamma_2}\subset U_{Q_0,Z_0}$ called the {\it walls corresponding to $Q_0,\gamma_1,\gamma_2$}. We denote their union by ${\cal W}_1:={\cal W}_1^{Q_0}$ and call it
{\it the wall of first kind} (physicists call it the {\it wall of marginal stability}).

Let us consider a continuous path $Z_t, 0\le t\le 1$ in $U_{Q_0,Z_0}$ which intersects each of these walls for finitely many values of $t\in [0,1]$. Suppose that we have a continuous lifting path $(Z_t,a_t)$ of this path in $Stab(\g)$ such that $Q_0$ is compatible with each $a_t$ for all $0\le t\le 1$. Then for any
$\gamma\in \Gamma\setminus \{0\}$ such that $Q_0(\gamma)\geqslant 0$ the element $a_t(\gamma)$ does not change as long as $t$ satisfies the condition $$Z_t(\gamma)\notin \cup_{\gamma_1,\gamma_2\in \Gamma\setminus \{0\},\, \gamma_1+\gamma_2=\gamma}{\cal W}_{\gamma_1,\gamma_2}^{Q_0}\,\,.$$ If this condition is not satisfied we say that $t$ is a discontinuity point for $\gamma$.
For a given $\gamma$ there are finitely many discontinuity points.

Notice that for each $t\in [0,1]$ there exist limits
$$a^{\pm}_t(\gamma)=\lim_{\varepsilon\to 0, \,\varepsilon>0}a_{t\pm \varepsilon}(\gamma)$$ (for $t=0$ or $t=1$ only one of the two limits is well-defined).
Then the continuity of the lifted path $(Z_t,a_t)$ is equivalent to the following {\it wall-crossing formula}
which holds for {\it any $t\in [0,1]$} and arbitrary $\gamma\in \Gamma\setminus\{0\}$:
$$\prod_{\mu\in \Gamma^{prim},\,Z_t(\mu)\in l_{{\gamma},t}}^{\longrightarrow}\exp\left(\sum_{n\geqslant 1}a_t^{-}(n\mu)\right)=$$
$${}=\exp\left(\sum_{\mu\in \Gamma^{prim},\,Z_t(\mu)\in l_{{\gamma},t}, \,n\geqslant 1}a_t(n\mu)\right)=\prod_{\mu\in \Gamma^{prim},\,Z_t(\mu)\in l_{{\gamma},t}}^{\longrightarrow}\exp\left(\sum_{n\geqslant 1}a^{+}_t(n\mu)\right),$$
where $l_{\gamma,t}=\R_{>0}Z_t(\gamma)$, and $\Gamma^{prim}\subset \Gamma$ is the set of primitive vectors.
The first and the last products are taken in the clockwise order of $\op{Arg}(Z_{t-\varepsilon})$ and
$\op{Arg}(Z_{t+\varepsilon})$ respectively, where $\varepsilon>0$ is sufficiently small.
Moreover, for each $\gamma$ we have $a_t^{-}(\gamma)=a_t^{+}(\gamma)=a_t(\gamma)$ unless there exist non-zero $\gamma_1,\gamma_2$ such that $\gamma=\gamma_1+\gamma_2$ and $Z_t\in {\cal W}_{\gamma_1,\gamma_2}^{Q_0}$.

\begin{rmk} Informally speaking, the wall-crossing formula says that for a very small sector $V$ containing the ray $l_{{\gamma},t}$ the corresponding element $A_V$, considered as a function of time, is locally constant in a neighborhood of $t$.

\end{rmk}

For each $\gamma\in \Gamma\setminus \{0\}$ the wall-crossing formula allows us to calculate $a_{1}(\gamma)$ is terms of $a_{0}(\gamma^{\prime})$ for a finite collection of elements $\gamma^{\prime}\in \Gamma\setminus \{0\}$. Morally it is an inductive procedure on the ordered set of discontinuity points $t_i\in [0,1]$.

In the case of a $3CY$ category the wall-crossing formula can be written in a very explicit form. For that it is not necessary to enlarge the ``charge lattice" $\Gamma$. More precisely,
let $\Gamma$ be a free abelian group of finite rank $n$, endowed
with a skew-symmetric integer-valued  bilinear form $\langle \bullet, \bullet\rangle : \Gamma \times \Gamma\to \Z$. Then  the corresponding Lie algebra of functions on a Poisson torus (which is the quasi-classical limit of the quantum torus associated with the $3CY$ category) is given explicitly in terms of generators and relations as
$\g_{\Gamma}=\g_{\Gamma,\langle \bullet,\bullet \rangle}=\oplus_{\gamma\in \Gamma}\Q\cdot e_{\gamma}$ with the Lie bracket
$$[e_{\gamma_1},e_{\gamma_2}]=(-1)^{\langle \gamma_1, \gamma_2\rangle}\langle \gamma_1, \gamma_2\rangle e_{\gamma_1+\gamma_2}\,\,.$$

It can be made into a Poisson algebra with the commutative prouct given by
$e_{\gamma_1}e_{\gamma_2}=(-1)^{\langle \gamma_1,\gamma_2\rangle}\,e_{\gamma_1+\gamma_2}.$

For a stability data $(Z,a)$ we can write uniquely (by the M\"obius inversion formula)
$$a(\gamma)=-\sum_{n\geqslant 1, {1\over{n}}\gamma\in \Gamma\setminus \{0\}}{\Omega(\gamma/n)\over{n^2}}e_{\gamma}\,\,,$$
where $\Omega:\Gamma\setminus \{0\}\to \Q$ is a function. Then we have
$$\exp\left(\sum_{n\geqslant 1}a(n\gamma)\right)=\exp\left(-\sum_{n\geqslant 1} \Omega(n\gamma)\sum_{k\geqslant 1}{e_{kn\gamma}\over{k^2}}\right):=
\exp\left(-\sum_{n\geqslant 1}\Omega(n\gamma)\op{Li}_2(e_{n\gamma})\right),$$
where $\op{Li}_2(t)=\sum_{k\geqslant 1}{t^k\over{k^2}}$ is the dilogarithm function.
Then
$$ A_V:=\prod_{\gamma\in C(V,Z)\cap \Gamma}^{\longrightarrow}
    \exp\left(- \Omega(\gamma)\sum_{n=1}^\infty \frac{e_{n\gamma}}{n^2}\right).$$

If the stability data belongs to the wall of first kind then there exists a $2$-dimensional lattice $\Gamma_0\subset \Gamma$  such that $Z_{t_0}(\Gamma_0)$ belongs to a real line ${\R}e^{i\alpha}$ for some $\alpha\in [0,\pi]$.

The wall-crossing formula describes  the change of the values $\Omega(\gamma)$ for $\gamma \in \Gamma_0$ and depends only on the restriction
 $\Omega_{|\Gamma_0}$ of $\Omega$ to the  lattice $\Gamma_0$. Values $\Omega(\gamma)$ for $\gamma\notin \Gamma_0$
do not change at $t=t_0$.
Denote by $k\in \Z$ the value of the form $\langle\bullet,\bullet\rangle$ on a fixed basis $\gamma_1,\gamma_2$ of $\Gamma_0\simeq \Z^2$ such that $C(V)\cap \Gamma_0\subset \Z_{\geqslant 0}\cdot\gamma_1\oplus \Z_{\geqslant 0}\cdot\gamma_2$, where $C(V)$ is the convex cone defined previosuly. We assume that $k\ne 0$, otherwise there will be no jump in values of $\Omega$ on $\Gamma_0$.
The group  elements which we are interested in
can be identified with products of the following
automorphisms of $\Q[[x,y]]$ preserving the symplectic form $k^{-1}(xy)^{-1} dx\wedge dy$:
$$
T_{a,b}^{(k)}:(x,y)\mapsto $$
$$\mapsto \left( x\cdot (1-(-1)^{kab} x^a y^b)^{-kb}, y\cdot(1-(-1)^{kab} x^a y^b)^{ka}\right),\medskip a,b\geqslant 0, a+b\geqslant 1\,.
$$
For $\gamma=a\gamma_1+b\gamma_2$ we have
$$T_{a,b}^{(k)}=\exp\left(-\sum_{n\geqslant 1}\frac{e_{n\gamma}}{n^2}\right)$$
in the above notation.
Any exact symplectomorphism $\phi$ of $\Q[[x,y]]$ can be decomposed uniquely into a clockwise and an anti-clockwise product:
$$\phi=\prod_{a,b}^{\longrightarrow} \left(T_{a,b}^{(k)}\right)^{c_{a,b}}=\prod_{a,b}^{\longleftarrow} \left(T_{a,b}^{(k)}\right)^{{d}_{a,b}}$$
with certain exponents $c_{a,b},d_{a,b}\in \Q$. These exponents should be interpreted as the
limiting values
of the functions $\Omega_{t_0}^{\pm}=\lim_{t\to t_0\pm 0}\Omega_t$ restricted to $\Gamma_0$.
The passage from the clockwise order (when the slope $a/b \in [0,+\infty]\cap {\mathbb P}^1(\Q)$ decreases) to the anti-clockwise order (when the slope increases) gives the change of $\Omega_{|\Gamma_0}$ as we cross the wall.
It will be convenient to denote $T_{a,b}^{(1)}$ simply by $T_{a,b}$. The pronilpotent group generated by
the transformations $T_{a,b}^{(k)}$ coincides with the one generated by the transformations $T_{a,|k|b}$.

For instance the decompositions for $k=1,2$ look such as follows:
$$T_{1,0} \cdot T_{0,1} =T_{0,1}\cdot T_{1,1}\cdot T_{1,0}\,\,,$$
$$T_{1,0}^{(2)} \cdot T_{0,1}^{(2)}= T_{0,1}^{(2)}\cdot T_{1,2}^{(2)}\cdot T_{2,3}^{(2)}\cdot\dots
    \cdot (T_{1,1}^{(2)})^{-2} \cdot \dots\cdot T_{3,2}^{(2)} \cdot T_{2,1}^{(2)}
\cdot T_{1,0}^{(2)}\,\,,$$
or equivalently
$$T_{1,0} \cdot T_{0,2}= T_{0,2}\cdot T_{1,4}\cdot T_{2,6}\cdot\dots
    \cdot T_{1,2}^{-2} \cdot \dots\cdot T_{3,4} \cdot T_{2,2}
\cdot T_{1,0}\,\,.$$

The formula for $k=2$ is related to the one for the BPS spectrum in the pure $N=2, d=4$ Seiberg-Witten theory.

\subsection{Walls of second kind}

Let us fix as above the quadratic form $Q_0$ and a connected component $U$ of the set $\{Z\in \op{Hom}(\Gamma,\C)|\,\, (Q_{0})_{|\op{Ker} Z}<0\}$.
For a given primitive $\gamma\in \Gamma\setminus \{0\}$ we introduce the set ${\cal W}_{\gamma}^{Q_0}=\{Z\in U|\,Z(\gamma)\in \R_{>0}\}$. It is a real hypersurface in $U$. We call it the {\it wall of second kind} associated with $\gamma$. We call the union $\cup_{\gamma}{\cal W}_{\gamma}^{Q_0}$ the {\it wall of second kind} and denote it by
${\cal W}_2$.

\begin{defn} We say that a path $\sigma=(Z_t)_{0\le t\le 1}\subset U$ is short if the convex cone
$C_{\sigma}$ which is the convex hull of $\left(\cup_{0\le t\le 1}Z_t^{-1}(\R_{>0})\right)\cap \{Q_0\geqslant 0\}$
is strict.

\end{defn}

With a short path we  associate a pronilpotent group $G_{C_{\sigma}}$ with  Lie algebra
$\g_{C_{\sigma}}=\prod_{\gamma\in C_{\sigma}\cap \Gamma}\g_{\gamma}$.

It follows from the definition of the topology that there is a lifting map
$\phi: U\to Stab(\g)$ (more precisely, one uses the fact that the projection of the space of stability conditions to the space of central charges is a local homeomorphism).

Then one can easily prove the following result.

\begin{prp} For a generic short path $\sigma=(Z_t)_{0\le t\le 1}$ there exists a no more than countable
set $t_i\in [0,1]$ and corresponding primitive $\gamma_i\in \Gamma\setminus \{0\}$ such that $Z_{t_i}\in {\cal W}_{\gamma_i}^{Q_0}$. For each $i$ we have: $\op{rk} Z_{t_i}^{-1}(\R)\cap \Gamma=1$.

\end{prp}

For every such $t_i$ we introduce
$$A_{t_i}=\exp\left(\varepsilon_i\sum_{n\geqslant 1}a_{t_i}({n\gamma_i})\right)\in G_{C_{\sigma}}\,\,,$$
where $\varepsilon_i=\pm 1$ depending on the direction in which the path $Z_t(\gamma_i)$ crosses $\R_{>0}$ for $t$ sufficiently close to $t_i$.

Then the wall-crossing formulas  imply the following result (Theorem 4 from Section 2.4 of [KoSo1]).

\begin{thm} For any short loop the monodromy $\prod_{t_i}^{\longrightarrow}A_{t_i}$ is equal to the identity
(here the product is taken in the increasing order of the elements $t_i$).

\end{thm}

This ``triviality of the monodromy"  allows us to introduce (at least for a short path) the well-defined notion of parallel transport. In the categorical framework the parallel transport appears when one calculates the change of the ``numerical" DT-invariants under the change of $t$-structure (e.g. under tilting).

\section{Examples}

\subsection{Quivers with potential and cluster transformations}

Let ${\cal C}$ be a $3$-dimensional  ind-constructible Calabi-Yau category over a field $\k$ of characteristic zero. Assume that it is endowed with a finite collection of spherical generators
${\cal E}=\{E_i\}_{i\in I}$ of ${\cal C}$ defined over $\k$.  Then
$\op{Ext}_{\cal C(\k)}^\bullet(E_i,E_i)$ is
isomorphic to $H^{\bullet}(S^3,\k ), \,\,i\in I$.
The  matrix of the Euler form (taken with the {\it minus} sign)
$$a_{ij}:=-\chi\left(\op{Ext}_{\cal C(\k)}^{\bullet}(E_i,E_j)\right)$$
is integer and skew-symmetric.
 In fact, the ind-constructible category   ${\cal C}$ can be canonically derived from the
$\k$-linear Calabi-Yau $\A$-category $\CC(\K)$, or even from its full subcategory consisting of the collection ${\cal E}$.
 In what follows we will omit the subscript $\CC(\k)$ in the notation for $\op{Ext}^\bullet$-spaces.

\begin{defn} The collection ${\cal E}$ is called a cluster if for any $i\ne j$ the graded space $\oplus_{m\in \Z} \op{Ext}^m(E_i,E_j)$ is either zero, or it is concentrated
 in one of the two degrees $m=1$ or $m=2$ only.

\end{defn}

We will assume that our collection is a cluster. In that case $K_0(\CC(\kk))\simeq \Z^I$ with the basis formed by the isomorphism classes $[E_i], i\in I$. With a cluster collection one associates a quiver $Q$ with potential $W$ in the natural way: the vertices of $Q$ correspond to the cluster generators, the arrows between vertices $i$ and $j$ correspond to a basis in $\op{Ext}^1(E_i,E_j)$ and the potential $W$ is the restriction of the potentail of the category to the $\oplus_{i,j\in I} \op{Ext}^1(E_i,E_j)$. It is proved in [KoSo1], Section 8.1, Theorem 9 that the pairs $(Q,W)$ up to a naturally defined gauge equivalence correspond to  $3CY$ categories (up to equivalence) generated by a cluster collection for which
\begin{itemize}
\item
 $\op{Ext}^0(E_i,E_i)=\k \,id_{E_i}$,
\item $\op{Ext}^0(E_i,E_j)=0$ for any $i\ne j$,
\item $\op{Ext}^{<0}(E_i,E_j)=0$, for any $i,j$.
\end{itemize}

Mutation (a.k.a. spherical reflection) in $\CC$ corresponds to the notion of mutation of $(Q,W)$ in quiver theory. Cluster generators single out an open domain in $Stab(\CC)$, where all central charges $Z(E_i)$ belong to the upper half-plane. Categorical mutations change the $t$-structure, so that some $Z(E_i)$ can move to the lower half-plane. It can be geometrically interpreted as a rotation of the upper-half plane and hence corresponds to the crossing of the wall of second kind. The motivic DT-invariant gives rise to an automorphism of the corresponding quantum torus, as we discussed before. It can be written in terms of the quiver $(Q,W)$ as well as in terms of the mutated quiver $(Q^{\prime}, W^{\prime})$. Two expressions of the same motivic DT-invariant are related by a quantum cluster transformation (see [KoSo1], Section 8.4 Corollary 4 for a precise statement). It would be interesting to compare this picture with the one developed in [GaMoNe2].

\subsection{Stability data from complex integrable systems}

In Section 2.7 of [KoSo1] we explained how one can associate stability data on a graded Lie algebra to a complex integrable system.
Recall that a {\it complex integrable system} is a holomorphic map $\pi:X\to B$ where $(X,\omega^{2,0}_X)$ is a holomorphic symplectic manifold, $\dim X=2\dim B$,  and the generic fiber of $\pi$ is a Lagrangian submanifold, which is a polarized abelian variety. We assume (in order to simplify the exposition) that the polarization is
principal. The fibration $\pi$ is non-singular outside of a closed subvariety $B^{sing}\subset B$ of complex codimension at least one.
It follows that on the open subset $B^{sm}:=B\setminus B^{sing}$ we have a local system ${\bf \Gamma}$ of symplectic lattices with the fiber over $b\in B^{sm}$ equal to $\Gamma_b:=H_1(X_b,\Z), X_b=\pi^{-1}(b)$ (the symplectic structure on $\Gamma_b$ is given by the  polarization).

Furthermore, the set $B^{sm}$ is locally (near each point $b\in B^{sm}$) embedded as a holomorphic Lagrangian subvariety into an affine symplectic space parallel to $H_1(X_b,\C)$. Namely, let us choose a symplectic basis $\gamma_i\in \Gamma_b, 1\le i\le 2n$. Then we have a collection of holomorphic closed
$1$-forms $\alpha_i=\int_{\gamma_i}\omega^{2,0}_X, 1\le i\le 2n$ in a neighborhood of $b$. There exists (well-defined locally up to an additive constant) holomorphic functions $z_i, 1\le i\le 2n$ such that $\alpha_i=dz_i, 1\le i\le 2n$. They define an embedding of a neighborhood of $b$ into $\C^{2n}$. The collection of $1$-forms $\alpha_i$ gives rise to an element $\delta\in H^1(B^{sm}, {\bf \Gamma}^{\vee}\otimes \C)$. {\it We  assume that $\delta=0$}.
This assumption is equivalent to an existence of a section
$Z\in \Gamma(B^{sm},{\bf \Gamma}^{\vee}\otimes {\cal O}_{B^{sm}})$ such that $\alpha_i=dZ(\gamma_i), 1\le i\le 2n$.
\begin{defn} We call $Z$ the central charge of the integrable system.

\end{defn}

Hence, for every point $b\in B^{sm}$ we have a symplectic lattice $\Gamma_b$ endowed with an additive map
$Z_b:\Gamma_b\to \C$. In the loc. cit. we defined a continuous family of stability data on graded Lie algebras
$\mathfrak{g}_{\Gamma_b}$ with central charges $Z_b$ in the following way.

First, observe that the dense open set $B^{sm}\subset B$ carries a K\"ahler form
$$\omega_{B}^{1,1}=\op{Im}\left(\sum_{1\le i\le n}\alpha_i\wedge\overline{\alpha_{n+i}}\right)\,\,.$$ We denote by $g_{B}$ the corresponding K\"ahler metric.

For any $t\in \C^{\ast}$ we define an integral affine structure on the $C^{\infty}$-manifold $B^{sm}$ given by a collection of closed $1$-forms $\op{Re}(t\alpha_i), 1\le i\le 2n$. For any simply-connected open subset $U\subset B^{sm}$ and a covariantly constant section $\gamma\in \Gamma(U,{\bf \Gamma})$ we have a closed $1$-form
$$\alpha_{\gamma,t}=\op{Re}\left(t\int_{\gamma}\omega_X^{2,0}\right)=d\op{Re}(tZ(\gamma))\,\,,$$
and the corresponding gradient vector field $v_{\gamma,t}=g_B^{-1}(\alpha_{\gamma,t})$.
Notice that this vector field is a constant field with integral direction in the integral affine structure associated with the closed $1$-forms $\op{Im}(t\alpha_i), 1\le i\le 2n$.
Second, using the approach of [KoSo4] we can construct an infinite oriented tree lying in $B$ such that its external vertices belong to $B^{sing}$, and its edges are {\it positively oriented } trajectories of the vector fields $v_{\gamma,t}$. All internal vertices have valency at least $3$, and  every such vertex should be thought of as a splitting point: a trajectory of the vector field $v_{\gamma,t}$ is split at a vertex into several trajectories of some vector fields $v_{\gamma_1,t},\dots,v_{\gamma_k,t}$ such that $\gamma=\gamma_1+\dots+\gamma_k$.

For a fixed $t\in \C^{\ast}$ the union $W_t$ of all trees as above is in fact a countable union of real hypersurfaces in $B^{sm}$. They are analogs of the walls of second kind. The set $W_t$ depends on $\op{Arg}t$ only. The union
$\cup_{\theta\in [0,2\pi )}W_{te^{i\theta}}$ swaps the whole space $B^{sm}$. Let us denote by $W^{(1)}$ the union over all $t\in \C^{\ast}/\R_{>0}$ of the sets of internal vertices of all trees in $W_t$ (splitting points of the gradient trajectories). This is an analog of the wall of first kind. Once again using the approach of [KoSo4] we assign rational multiplicities to edges of the tree.
 This leads to the following picture. Consider the total space $tot({\bf \Gamma})$ of the local system ${\bf \Gamma}$. It follows from above assumptions and considerations that we have a locally constant function $\Omega:tot({\bf \Gamma})\to \Q$ which jumps at the subset consisting of the lifts of the wall $W^{(1)}$ to $tot({\bf \Gamma})$.  Then for a fixed $b\in
B^{sm}$ the pair $(Z,\Omega)$ defines stability data on the graded Lie algebra $\g_{\Gamma_b}$ of the group of formal symplectomorphisms of the symplectic torus associated with $\Gamma_b$. In this way we obtain a local embedding $B^{sm}\hookrightarrow Stab(\g_{\Gamma_b})$.

Examples of this construction include integrable systems of Seiberg-Witten theory, Hitchin systems, etc. The case of pure $SU(2)$ Seiberg-Witten theory was illustrated in Section 2.7 of [KoSo1]. The corresponding wall-crossing formula
coincides with the one for $T_{a,b}^{(2)}$ considered above.
More examples can be found in [GaMoNe1,2]. We also remark that this section is related to  Sections 7.2, 7.6 below.

\subsection{Cecotti-Vafa work and WCF for $gl(n)$}

This is an example of the wall-crossing formula considered in [KoSo1], Section 2.9.

Let $\g={\mathfrak{gl}}(n,\Q)$ be the Lie algebra of the general linear group. We consider it as a $\Gamma$-graded Lie algebra $\g=\oplus_{\gamma \in \Gamma}\g_{\gamma}$,  where $$\Gamma=\{(k_1,\dots,k_n)|\,\,k_i\in \Z, \sum_{1\le i\le n}k_i=0\}$$ is the root lattice. We endow $\g$ with the Cartan involution $\eta$. The algebra $\g$ has
 the standard basis $E_{ij}\in \g_{\gamma_{ij}}$  consisting of matrices with the single non-zero entry at the place $(i,j)$ equal to $1$. Then $\eta(E_{ij})=-E_{ji}$. In [KoSo1], Section 2.1 we introduced also the notion of symmetric stability data on a graded Lie algebra. The definition includes a choice of involution $\eta$ on $\Gamma$. In this section  we are going to consider symmetric stability data.

We notice that  $$\op{Hom}(\Gamma,\C)\simeq \C^n/\C\cdot(1,\dots,1)\,\,.$$ We define a subspace $\op{Hom}^\circ(\Gamma,\C)\subset \op{Hom}(\Gamma,\C)$ consisting (up to a shift by the multiples of the vector $(1,\dots,1)$) of vectors $(z_1,\dots,z_n)$ such that $z_i\ne z_j$ if $i\ne j$. Similarly we define a subspace $\op{Hom}^{\circ\circ}(\Gamma,\C)\subset \op{Hom}(\Gamma,\C)$ consisting (up to the same shift) of such $(z_1,\dots,z_n)$ that there is no $z_i,z_j,z_k$ belonging to the same real line as long as $i\ne j\ne k$. Obviously there is an inclusion $\op{Hom}^{\circ\circ}(\Gamma,\C)\subset \op{Hom}^\circ(\Gamma,\C)$.

For $Z\in \op{Hom}(\Gamma,\C)$ we have $Z(\gamma_{ij})=z_i-z_j$. If $Z\in \op{Hom}^{\circ\circ}(\Gamma,\C)$ then symmetric stability data with such $Z$ is the same as a skew-symmetric matrix $(a_{ij})$ with rational entries determined from the equality $a(\gamma_{ij})=a_{ij}E_{ij}$. Every continuous path in $\op{Hom}^\circ(\Gamma,\C)$ admits a unique lifting to $Stab(\g)$ as long as we fix the lifting of the initial point. The matrix $(a_{ij})$ changes when we cross walls in  $\op{Hom}^\circ(\Gamma,\C)\setminus \op{Hom}^{\circ\circ}(\Gamma,\C)$. A typical wall-crossing corresponds to the case when in the above notation the point $z_j$ crosses a straight segment joining $z_i$ and $z_k, i\ne j\ne k$. In this case the only change in the matrix $(a_{ij})$ is of the form:
$$a_{ik}\mapsto a_{ik}+a_{ij}a_{jk}\,\,.$$
Exactly these wall-crossing formulas appeared in [CeVa] (the numbers $a_{ij}$ are integers in the loc. cit.)

\subsection{About explicit formulas}

Assume that a $3$-dimensional Calabi-Yau category $\CC$ is generated by one spherical object $E$ defined over $\k$. Therefore $R:=\op{Ext}^{\bullet}(E,E)\simeq H^{\bullet}(S^3,\k)$. In this case we take $\Gamma=K_0(\CC(\kk))\simeq \Z$, and the skew-symmetric form on $\Gamma$ is trivial. We explained the choice of orinetation data in [KoSo1], Section 6.4.
For any $z\in \C, \op{Im}z>0$ we have a stability condition $\sigma_z$ such that $E\in \CC^{ss}, \,Z(E):=Z(\op{cl}_\kk(E))=z,\, \op{Arg}(E)=\op{Arg}(z)\in (0,\pi)$.
For a strict sector $V$ such that $\op{Arg}(V)\subset (0,\pi)$ we have the category $\CC_V$ which is either trivial (if $z\notin V$) or consists of objects $0,E,E\oplus E,\dots$ (if $z\in V$).
Then $A_V^{mot}=1$ in the first case and
$$A_V^{mot}=\sum_{n\geqslant 0}\frac{{\mathbb L}^{n^2/2}}{[GL(n)]}\hat{e}_{\gamma_1}^n\,\,,$$
in the second case. Here  $\gamma_1$ is the generator of $\Gamma$ (i.e. the class of $E$ in the $K$-theory).

In this case $\op{Ext}^1(nE,nE)=0$, where we set $nE=E^{\oplus n}, n\geqslant 1$. Therefore  $W_{nE}=0$ which implies that for the motivic Milnor fiber we have $MF(W_{nE})=0$.
The numerator in the above formula is
$${\mathbb L}^{n^2/2}={\mathbb L}^{{1\over{2}}\,\dim \op{Ext}^0(nE,nE)}=
{\mathbb L}^{{1\over{2}}\sum_{i\le 1}(-1)^i \dim \op{Ext}^i(nE,nE)}\,\,,$$
since $\op{Ext}^{1}(nE,nE)=0, \op{Ext}^{<0}(nE,nE)=0$.

Let us consider the ``quantum dilogarithm" series
$${\bf E}(q^{1/2},x)=\sum_{n\geqslant 0}{{q}^{n^2/2}\over{(q^n-1)\dots(q^n-q^{n-1})}}x^n\in \Q(q^{1/2})[[x]]\,\,.$$

Since $[GL(n)]=({\mathbb L}^n-1)\dots({\mathbb L}^n-{\mathbb L}^{n-1})$, we conclude that
$$A_V^{mot}={\bf E}({\mathbb L}^{1/2},\hat{e}_{\gamma_1})\,\,.$$

In Section 6.4 of [KoSo1] we presented wall-crossing formulas for both motivic and numerical DT-invariants in  several cases (e.g. for D0-D6 BPS bound states). There are many more computations of this sort in the literature. Most of them are related to toric local $3CY$ varieties and the count often reduces to the count of cyclic modules similarly to [KoSo1], Sections 7.3,7.4. We mention here only few papers which contain explicit formulas: [DeMo], [GuDi], [Jaf], [JafMo], [NagNak],[Nag1], [Nag2], [OoYa], [Sz].

\section{Discussion, speculations, open questions}

There are several foundational questions in the theory of motivic DT-invariants which have to be settled, e.g. deformation invariance, existence of orientation data (or its invariance under mutations in the case of quivers with potential), existence of the quasi-classical limit, proof of the integral identity in the general motivic case. All those questions were addressed in [KoSo1]. We would like to recall below several other open problems discussed in [KoSo1] which, in our opinion, deserve further study.

\subsection{Intermediate Jacobian, complex integrable systems for 3CY and GW=DT}

The idea of associating a complex integrable system to a homologically smooth $3CY$ category was discussed in [KoSo1], Section 7.2.
It is partially motivated by the corresponding pure geometric story (see [DonMar]). In addition to that we speculated on how Donaldson-Thomas (and Gromov-Witten) invariants could emerge in this case.

More precisely, we suggested that for an arbitrary triangulated compact homologically smooth $\A$-category $\CC$ (see [KoSo2]) one has a non-commutative version of the Deligne
cohomology $H_{D}(\CC)$ which  fits into a short exact sequence defined via the Hodge filtration on periodic cyclic homology $HP_{\bullet}(\CC)$ (we also use the hypothetical topological $K$-group $K^{top}_{\bullet}$):
$$0\to HP_{odd}(\CC)/(F^{1/2}_{odd}+ K^{top}_{odd}(\CC))\to H_D(\CC)\to F^0_{even}\cap K^{top}_{even}(\CC)\to 0\,\,.$$
Morally, $H_D(\CC)$ should be thought of as zero cohomology group of the homotopy colimit of the following diagram of cohomology theories:
$$\begin{CD}  @. HC_\bullet^-(\CC)\\
  @. @VVV  \\
K^{top}_\bullet(\CC) @>>>  HP_\bullet(\CC)
\end{CD} $$
where $HC_\bullet^-(\CC)$ is the negative cyclic homology.

Any object of $\CC$ should have its characteristic class in $H_{D}(\CC)$. More precisely, there should be a homomorphism of groups $ch_D: K_0(\CC)\to H_{D}(\CC)$ (in the case of  Calabi-Yau manifolds it is related to holomorphic Chern-Simons functional). The reason for this is that every object $E\in Ob(\CC)$
should have  natural characteristic classes in $K^{top}_0(\CC)$ and in $HC_0^-(\CC)$ whose images in $HP_0(\CC)$ coincide with each other.

Let now $\CC$ be a $3CY$ category. Then we have the moduli space ${\cal M}$ of its deformations (including the Calabi-Yau structure). We conjectured in [KoSo1] that there is a fibration
${\cal M}^{tot}\to {\cal M}$
with the fiber $H_{D}(\CC)$ over the point $[\CC]\in {\cal M}$ such that ${\cal M}^{tot}$ is a holomorphic symplectic manifold. Moreover, any fiber of this fibration (i.e. the group $H_{D}(\CC)$ for given $[\CC]$) is a countable union of complex Lagrangian tori.
By analogy with the commutative case we suggested that the locus ${\cal L}\subset {\cal M}^{tot}$ consisting of values of $ch_D$ is a countable union of Lagrangian subvarieties. Every such subvariety can be either a finite ramified covering of ${\cal M}$ or a fibration over a proper subvariety of ${\cal M}$ with the fibers which are abelian varieties.

For generic $[\CC]\in {\cal M} $ one can use the triple $(K_0(\CC),H_{D}(\CC),ch_D)$ as a  triple $(K_0(\CC),\Gamma,\op{cl})$. Analogs of our motivic Donaldson-Thomas invariants $A_V^{mot}\in {\cal R}_V$ will be  formal countable sums of points in $H_{D}(\CC)$ with ``weights" which are elements of the motivic ring ${\overline{D^{\mu}}}$. The pushforward map from
$H_{D}(\CC)$ to $\Gamma=F^0_{even}\cap K^{top}_{0}(\CC)$ gives the numerical DT-invariants. The continuity of motivic DT-invariants means that after taking the quasi-classical limit the weights become integer-valued functions on the set of those irreducible components of ${\cal L}$ which are finite ramified coverings on ${\cal M}$.

These considerations lead to the following question raised in [KoSo1], Section 7.2

\begin{que} Is there a natural extension of the numerical DT-invariants to those components of ${\cal L}$ which project to a proper subvariety of ${\cal M}$?
\end{que}

The theory of Gromov-Witten invariants can be (again hypothetically) expressed in a similar way. This leads to a natural question about the generalization of the famous conjecture ``GW=DT".
The following remark from the loc. cit. can be considered as a proposal of such generalization in the geometric case.

\begin{rmk}  Suppose $X$ is a $3d$ complex compact Calabi-Yau manifold with $H^1(X,\Z)=0$. Then we have an exact sequence
$$0\to H^3_{DR}(X)/(F^2H^3_{DR}(X)+H^3(X,\Z))\to H^4_D(X)\to H^4(X,\Z)\to 0\,\,,$$
where $H^4_D(X)={\mathbb H}^4(X,\Z\to {\cal O}_X\to \Omega^1_X)$ is the Deligne cohomology.
Then any curve $C\subset X$ defines the class $[C]\in H^4_D(X)$. For a generic complex structure on $X$ the class is constant in any smooth connected family of curves. Moreover, a stable map to $X$ defines a class in $H^4_D(X)$. Then we have exactly the same picture with holomorphic symplectic fibration ${\cal M}^{tot}\to {\cal M}$ and Lagrangian fibers, as we discussed above. Similarly to the case of DT-invariants the GW-invariants appear as infinite linear combinations of points in $H^4_D(X)$, but this time with rational coefficients. We expect that the well-known relationship ``GW=DT" should be a statement about the equality of the above-discussed counting functions (assuming positive answer to the Question 1).

\end{rmk}

\subsection{Non-archimedean integrable systems}

In [KoSo1], Section 1.5 we  suggested that the moduli space of stability conditions (or, rather a Lagrangian cone in it with the induced affine structure) should be thought of as a base of non-archimedean analytic integrable system similar to the one considered in [KoSo4].
Indeed, the collection of formal symplectomorphisms $A_V$ encoding the numerical DT-invariants give rise to a rigid analytic space ${\cal X}^{an}$ over any non-archimedean field, similarly to [KoSo4]. This space carries an analytic symplectic form and describes ``the behavior at infinity" of a (possibly non-algebraic) formal smooth symplectic scheme over $\Z$. Motivated by the String Theory we conjectured that there exists an actual complex symplectic manifold ${\cal M}$ (vector or hyper multiplet moduli space) admitting a (partial) compactification
$\overline{\cal M}$ and such that its $\C((t))$-points are given by
$${\cal X}^{an}(\C((t)))=\overline{\cal M}(\C[[t]])\setminus ({\cal M}(\C[[t]])\cup (\overline{\cal M}\setminus {\cal M})(\C[[t]]))\,\,,$$
i.e. it is the space of formal paths hitting the compactifying divisor but not belonging to it).
In the case when $\CC$ is the Fukaya category of a complex $3d$ Calabi-Yau manifold $X$ the space ${\cal M}$ looks ``at infinity" as a deformation of a complex symplectic manifold ${\cal M}^{cl}$ where $\dim  {\cal M}=\dim  {\cal M}^{cl}=\dim  H^3(X,\C)$. The latter is the total space of the bundle ${\cal M}^{cl}\to {\cal M}_X$, where
${\cal M}_X$ is the moduli space of complex structures on $X$. The fiber of the bundle is isomorphic to the space
$$(H^{3,0}(X)\setminus \{0\})\times (H^3(X,\C)/H^{3,0}(X)\oplus H^{2,1}(X)\oplus H^3(X,\Z))$$
 parametrizing  pairs {\it (holomorphic volume element, point of the intermediate Jacobian)}.
As we explained above, we expect that there is a complex integrable system associated with an arbitrary homologically smooth $3d$ Calabi-Yau category. The fiber is the torus associated with the ``Deligne cohomology" of the category.
In some cases (see [GaMoNe1,2]) the total space carries more structures. We discuss some of them below in Section 7.6.
Furthermore, the total space (constructed in [KoSo4] as a non-archimedean analytic space) can be quantized (see [So1]). It would be interesting to reveal the physical meaning of this quantization at the level of hyperk\"ahler manifolds considered in [GaMoNe1,2].

\subsection{Motivic DT invariants and Cohomological Hall algebra}

An alternative approach to motivic DT-invariants was suggested in [KoSo3]. Technically speaking, it is presented in [KoSo3] in the case of quiver with potential only. But in fact the approach of [KoSo3] is more general, with the main formalism valid for so-called smooth representation towers introduced in the loc. cit. The latter structure is probably hidden in all $3CY$ categories which appear ``in nature". Since it is not our purpose  to explain here the results of [KoSo3], we will only illustrate them in the case of a quiver without potential.

Let us fix a finite quiver $Q$, with the set $I$ of vertices, and $a_{ij}\in \Z_{\ge 0}$ arrows from $i$ to $j$
 for $i,j\in I$.
 For any dimension vector
 $$\gamma=(\gamma^i)_{i\in I} \in \Z_{\ge 0}^I$$
 we have the space of representations of $Q$ in complex coordinate vector spaces of dimensions $(\gamma^i)_{i\in I}$:
 $$\M_\gamma\simeq \prod_{i,j\in I} \C^{a_{ij}\gamma^i \gamma^j}\,\,.$$

 This space is naturally endowed with the conjugate action of the complex algebraic group
  $$\mathbb {G}_\gamma:=\prod_{i\in I} \op{GL}(\gamma^i,\mathbb{C}).$$

  We use the standard model $\op{Gr}(d,\C^\infty):=\underrightarrow{\lim}\op{Gr}(d,\C^N),\,\,N\to +\infty$ of
the classifying space of $\op{GL}(d,\C)$ for $d\ge 0$, and define
  $$\B{\mathbb G}_{\gamma}:=\prod_{i\in I} \op{BGL}(\gamma^i,\C)=\prod_{i\in I} \op{Gr}(\gamma^i,\C^\infty)\,\,.$$
  Let us consider the universal family over  $\rm{B}{\mathbb G}_{\gamma}$
  $$\M_\gamma^{\univ}:=\left(\mathrm{E}{\mathbb G}_{\gamma}\times \M_\gamma\right)/{{\mathbb G}_{\gamma}}\,\,, $$ where
$\rm{E}{\mathbb G}_{\gamma}\to \B{\mathbb G}_{\gamma}$ is the universal
  ${\mathbb G}_\gamma$-bundle.
  We introduce a $\Z_{\ge 0}^I$-graded abelian group
  $${\cal H}:=\oplus_{\gamma} {\cal H}_\gamma\,\,,$$
  where each component is defined as an equivariant cohomology
  $${\cal H}_{\gamma}:=H^\bullet_{{\mathbb G}_{\gamma}} (\M_\gamma):=
  H^\bullet(\M_\gamma^{\univ})=\oplus_{i\ge 0} H^i(\M_\gamma^{\univ})\,\,.$$

Using naturally defined correspondences for the spaces $\M_\gamma$ we defined in [KoSo3] an associative product on
${\cal H}$, making it into  {\it cohomological Hall algebra} (it should be called ``generalized algebra of BPS states" in physics, since it is a mathematical incarnation of the BPS algebra discussed in [HaMo]). Toric localization gives an explicit formula for the product, which makes ${\cal H}$ into a special case of Feigin-Odessky shuffle algebra.

Let $Q=Q_d$ be now a quiver with just one vertex and $d\ge 0$ loops. Then the product formula specializes to the following one:

\begin{multline*}
(f_1\cdot f_2)(x_1,\dots,x_{n+m}):=\\
\sum_{\substack{ i_1<\dots<i_n\\ j_1<\dots<j_m\\
\{i_1,\dots,i_n,j_1,\dots,j_m\}=\\
=\{1,\dots,n+m\}}} f_1(x_{i_1},\dots,x_{i_n})\,f_2(x_{j_1},\dots,x_{j_m})\,
\left(\prod_{k=1}^n\prod_{l=1}^m(x_{j_l}-x_{i_k})\right)^{d-1}
\end{multline*}
for symmetric polynomials, $f_1$  in $n$  variables, and  $f_2$ in $m$ variables. The product $f_1\cdot f_2$ is a symmetric polynomial
in $n+m$ variables.

One can introduce a double grading on this algebra, by declaring for a homogeneous symmetric polynomial of degree $K$ in $n$ variables
 to have bigrading $(n,2K+(1-d)n^2)$. Equivalently, one can shift the cohomological grading in $H^\bullet(\mathrm{BGL}(n,\C))$ by
$[(d-1)n^2]$.

For general $d$, the Hilbert-Poincar\'e series $P_d=P_d(z,q^{1/2})$ of the bigraded algebra ${\cal H}$ is the generalized $q$-exponential function (a.k.a. as quantum dilogarithm):
$$\sum_{n\ge 0, m\in \Z}\dim({\cal H}_{n,m})\,z^n q^{m/2}=\sum_{n\ge 0} \frac{q^{(1-d)n^2/2}}
{(1-q)\dots(1-q^n)} z^n\in \Z((q^{1/2}))[[z]]\,\,.$$

Thus we see that
$$A_l(x,q)=P_d(x,q^{-1})\,\,.$$
Here in the LHS we consider the input of the ray $l$ in the upper-half plane to the motivic DT-invariants of the $3CY$ category with the heart consisting of the finite-dimensional representations of $Q$. In the case of non-zero potential the heart consists of its critical points, and one should use the cohomology of the perverse sheaf of vanishing cycles associated with the potential (see [KoSo3] for details). The cohomological Hall algebra does not depend on a choice of a stability condition (more precisely, it depends on a $t$-structure, but not on a central charge). As we show in [KoSo3] such a choice gives rise to a PBW-type basis in ${\cal H}$. Moreover it  also gives rise  to a factorization of $A_V$ into the product of factors of the type $A_{\gamma_j}^{ss}$, where $\gamma_j$ is a non-zero primitive element of $\Z^{I}$, and each factor $A_{\gamma_j}^{ss}$ is an infinite sum of equivariant cohomology of the moduli space of semistable representations of dimensions $n\gamma_j, n\ge 0$. This is exactly the product decomposition of $A_V$ which we discussed above for  general $3CY$ categories.

\subsection{Donaldson 4d theory, Borcherds automorphic forms}

We suggested in [KoSo1] that our wall-crossing formulas should be related to those in the Donaldson theory of $4d$ manifolds with $b_2^+=1$ (cf. with the approach via virtual fundamental class in [Moch]) as well as with Borcherds hyperbolic Kac-Moody algebras and multiplicative automorphic forms.
Since Donaldson theory depends on a choice of the gauge group, the latter should be somehow encoded in our story.
In case of a complex projective surface $S$ the expected relationship should combine the well-known description of Seiberg-Witten theory via complex integrable systems (see [Don]) with the DT-theory developed in [KoSo1] in the case of the local $3CY$ given by the total space of the anticanonical bundle of $S$.
It is also not clear for us how to relate the wall-crossing formulas in Donaldson theory with the theory of stability data on graded Lie algebras recalled in Section 5. We think it is an interesting problem.

\subsection{WKB approximation in complex domain and WCF}

There is a striking similarity between our wall-crossing formulas and identities for the Stokes automorphisms in the theory of WKB asymptotics for a second order operator on ${\bf P}^1$ (see e.g. [DelDiPh], Section 3).  The latter story was highly motivated by the works of \'Ecalle, Pham, Voros and others (see e.g. [DelDiPh]). There is an underlying graded Lie algebra of ``alien derivatives" in the story, which is probably closely related to the motivic Galois group from [KoSo3]. Despite of a certain similarity with [KoSo1], Section 2.8 as well as with [BrTL] and some parts of [GaMoNe2] this relationship is not clear for  us.

\subsection{$3CY$ categories, black holes and the  split attractor flow}

Here we clarify  the discussion in [KoSo1], Section 1.5, 1)-3). More details will appear in [KoSo5]. Cf. also with Section 7.2.

Let ${\cal M}_{CFT}$ be the ``moduli space" of  unitary $N=2$ superconformal field theories. It is believed that ${\cal M}_{CFT}\simeq {\cal M}_A\times {\cal M}_B$, where for CFTs associated with a $3CY$ manifold $X$ the moduli space  ${\cal M}_A$ is the space fo complexified K\"ahler structures on $X$ while {\cal ${\cal M}_B$ is the moduli space of complex structures on $X$. Recall that superstring theory predicts an existence of a family of ``compactifying to $10$-dimensions" of superconformal field theories parametrized by the $4$-dimensional space-time endowed with Einstein metric. The latter has singularities at black holes. Using the time invariance we obtain a metric $g$ on $\R^3\setminus \{x_1,...,x_n\}$, where $x_i$ are positions of black holes.
This family can be interpreted as a harmonic map $h: \R^3\setminus \{x_1,...,x_n\}\to {\cal M}_{CFT}$.

Let us assume that our CFT is of geometric origin and comes from a $3CY$ manifold $X$. Assume that the K\"ahler component of $h$ is constant. Then according to [De1], [DeGrRa] (see aslo [DeMo]) the set of pairs $(h,g)$ is in one-to-one correspondence with the set of   maps

$$\phi:\R^3\setminus \{x_1,...,x_n\}\to {\cal M}_X$$
(here ${\cal M}_X$ is the moduli space of complex structures on $X$)
coming from the following ansatz.
Namely, the map $\phi$ is obtained by the projectivization of the map
$\hat{\phi}:\R^3\setminus \{x_1,...,x_n\}\to \Lambda$, where $\R^3\setminus \{x_1,...,x_n\}$ is endowed with the flat Euclidean metric and $\Lambda$ is the Lagrangian cone of the moduli space of deformations of $X$ endowed with a holomorphic volume form (it is locally embedded into $H^3(X,\C)$ via the period map). We endow $\Lambda$ with an integral affine structure via the local homeomorphism $Im:\Lambda\to H^3(X,\R), (\tau, \Omega^{3,0}_{\tau})\mapsto Im(\Omega^{3,0}_{\tau})$, where $\tau\in {\cal M}_X$ is a complex structure on $X$, and $\Omega^{3,0}_{\tau}$ is the corresponding holomorphic volume form.
Then the ansatz comes from harmonic maps $\hat{\phi}$ which are locally of the form
$Im\circ \hat{\phi}(x)=\sum_{1\le i\le n}{\gamma_i\over{|x-x_i|}}+v_{\infty}$, where $\gamma_i, 1\le i\le n$ are elements of the charge lattice $\Gamma=H^3(X,\Z)$ (their meanings are the charges of black holes) and $v_{\infty}$ is the boundary condition ``at infinity" (see [De1]).
This gives us $\phi$.

The image of $\phi$ is an ``amoeba-shaped" $3$-dimensional domain in ${\cal M}_X$.  Conjecturally, connected components of the moduli space of maps $\phi$ with given $v_{\infty}, \gamma_i, 1\le i\le n$ are in one-to-one correspondence with split attractor trees (see [De1]).
The edges of such a tree are the gradient trajectories of the function $|Z(\gamma)|=|\int_{\gamma}\Omega^{(3,0)}|^2/|\int_X\Omega^{(3,0)}\wedge \overline{\Omega^{(3,0)}}|$ .  Any edge is locally a projection of an affine line in $\Lambda$ with the slope $\gamma\in \Gamma=H^3(X,\Z)$. The condition ``at infinity", $v_{\infty}\in \Lambda$ is normalized so that $vol(v_{\infty})=1$. If the split attractor flow (lifted from  ${\cal M}_X$ to $\Lambda$) starting at $v_{\infty}$ in the direction $\gamma$ hits the wall of marginal stability so that $\gamma=\gamma_1+\gamma_2+..+\gamma_k$ then all $\gamma_1,...,\gamma_k$ belong (generically) to a two-dimensional plane.

One can argue (see more in [KoSo5]) that using our wall-crossing formulas it is possible to find all $\Omega(\gamma):=\Omega(b,\gamma), b\in \Lambda,\gamma\in \Gamma$ starting with a collection of integers $\Omega(b_{\gamma},\gamma)$ at the ``attractor points" $b_{\gamma}\in \Lambda$ given by the equation $Im\,b_{\gamma}=\gamma$. The points $\C^{\ast}b_{\gamma}\in {\cal M}_X$ are external vertices of the split attractor trees. The wall-crossing formulas are used at the internal vertices of the trees for the computation of $\Omega(b,\gamma)$. The numbers $\Omega(b_{\gamma},\gamma)$ can be arbitrary. We expect that they are determined uniquely by the geometry of certain quaternion K\"ahler manifold (hypermultiplet moduli space).

For a local $3CY$ manifold the quaternion K\"ahler structure reduces to a hyperk\"ahler one. In this case there is an alternative description of the numbers $\Omega(\gamma)$ as the virtual numbers of holomorphic discs with the boundary on Lagrangian tori (see [KoSo4], [GaMoNe1,2]).
This leads to a question about the relationship of semistable objects of the $\A$-category associated with our local Calabi-Yau manifold and holomorphic curves in the hyperk\"ahler manifold. Counting in both cases gives the same numbers.

More generally, one can speculate that for a non-compact quaternion K\"ahler manifold with certain behavior ``at infinity" one can construct a $3CY$ category whose objects are ``quaternion curves" in the twistor space.

\subsection{Stability conditions on curves and moduli of abelian differentials}

Geometry similar to the one on the space of stability conditions  appears in the theory of moduli spaces
of holomorphic abelian differentials (see e.g. [Zo]). The moduli space of abelian differentials
is a complex manifold, divided by real ``walls" of codimension one into pieces glued from
convex cones. It also carries a natural non-holomorphic action of the group $GL^+(2,\R)$. There is an analog of the central charge $Z$ in the story. It is given
by the integral of an abelian differential over a path between marked points in a complex curve.
We conjectured in [KoSo1] that the  moduli space of abelian differentials associated with
a complex curve with marked points, is isomorphic
to the moduli space of stability structures on the Fukaya category of this curve.
As a byproduct one can obtain an example of a non-connected moduli space of stability conditions
(cf. with [KoZo]).

\vspace{3mm}

{\bf References}
\vspace{3mm}

[Be1] K. Behrend, Donaldson-Thomas invariants via microlocal geometry
arXiv:math/0507523.
\vspace{2mm}

[BeFa] K. Behrend, B. Fantechi, Symmetric obstruction theories and Hilbert schemes of points on threefolds, arXiv:math/0512556.

\vspace{2mm}

[BiMi] E. Bierstone, P. Milman, Functoriality in resolution of singularities, ArXiv:math/0702375.

\vspace{2mm}

[BonVdB] A. Bondal, M. Van den Bergh,  Generators and representability of functors in commutative and noncommutative geometry, arXiv:math/0204218.

\vspace{2mm}
[Br1] T. Bridgeland, Stability conditions on triangulated categories,
arXiv:math/0212237.

\vspace{2mm}

[BrTL] T. Bridgeland, V. Toledano Laredo, Stability conditions and Stokes factors, arXiv:0801.3974.

[CeVa] S. Cecotti, C. Vafa, On Classification of N=2 Supersymmetric Theories, arXiv:hep-th/9211097.

\vspace{2mm}

[De1] F.Denef, Supergravity flows and D-brane stability, arXiv:hep-th/0005049.

\vspace{2mm}

[DeMo] F. Denef, G. Moore,  Split States, Entropy Enigmas, Holes and Halos, arXiv:hep-th/0702146.

\vspace{2mm}

[DeGrRa] F. Denef, B. Greene, M. Raugas, Split attractor flows and the spectrum of BPS D-branes on the Quintic, arXiv:hep-th/0101135.

\vspace{2mm}

[DelDiPh] E. Delabaere, H. Dillinger, F. Pham, Resurgence de Voros et p\'eridodes des courbes hyperelliptiques, An. Inst. Fourier, 43:1, 163-199, 1993.

\vspace{2mm}

[DenLo] J. Denef, F. Loeser, Geometry on arc spaces of algebraic varieties, arXiv:math/0006050.
\vspace{2mm}

[Don] R. Donagi, Seiberg-Witten integrable systems, arXiv:alg-geom/9705010.

\vspace{2mm}

[DonMar] R. Donagi, E. Markman, Cubics, integrable systems, and Calabi-Yau threefolds, arXiv:alg-geom/9408004.

\vspace{2mm}

[GaMoNe1] D. Gaiotto, G. Moore, A. Neitzke,  Four-dimensional wall-crossing via three-dimensional field theory, arXiv:0807.4723.

\vspace{2mm}
[GaMoNe2] D. Gaiotto, G. Moore, A. Neitzke,  Wall-crossing, Hitchin Systems, and the WKB Approximation, arXiv:0907.3987.

\vspace{2mm}

[GuDi] S. Gukov, T. Dimofte,  Refined, Motivic, and Quantum,  arXiv:0904.1420.

\vspace{2mm}

[HaMo] J. A. Harvey, G. Moore,  Algebras, BPS States, and Strings, arXiv:hep-th/9510182.

\vspace{2mm}

[Jaf] D. Jafferis,  Topological Quiver Matrix Models and Quantum Foam,  arXiv:0705.2250.
\vspace{2mm}

[JafMo] D. Jafferis, G. Moore, Wall crossing in local Calabi Yau manifolds, arXiv:0810.4909.
\vspace{2mm}

[Jo1] D. Joyce, Configurations in abelian categories. I. Basic properties and moduli stacks, arXiv:math/0312190.

\vspace{2mm}

[Jo2] D. Joyce,  Configurations in abelian categories. II. Moduli stacks, arXiv:math/0312192.
\vspace{2mm}

[Jo3] D. Joyce, Constructible functions on Artin stacks,  arXiv:math/0403305.

\vspace{2mm}

[Jo4] D. Joyce, Configurations in abelian categories. IV. Invariants and changing stability conditions,  arXiv:math/0410268.

[JoS1] D. Joyce, Y. Song, A theory of generalized Donaldson-Thomas invariants, to appear.

\vspace{2mm}

[JoS2] D. Joyce, Y. Song, A theory of generalized Donaldson-Thomas invariants, arXiv:0810.5645
\vspace{2mm}

[Ke1] B. Keller,  Introduction to A-infinity algebras and modules,  arXiv:math/9910179.

\vspace{2mm}

[KoSo1] M. Kontsevich, Y. Soibelman,  Stability structures, motivic Donaldson-Thomas invariants and cluster transformations,arXiv:0811.2435.

\vspace{2mm}

[KoSo2] M. Kontsevich, Y. Soibelman, Notes on A-infinity algebras, A-infinity categories and non-commutative geometry. I, arXiv:math/0606241.

\vspace{2mm}

[KoSo3] M. Kontsevich, Y. Soibelman, Cohomological Hall algebra and motivic Donaldson-Thomas invariants, to appear.

\vspace{2mm}

[KoSo4] M. Kontsevich, Y. Soibelman, Affine structures and non-archimedean analytic spaces, arXiv:math/0406564.

\vspace{2mm}
[KoSo5] M. Kontsevich, Y. Soibelman, Split attractor flow, BPS invariants and affine geometry, work in progress.

\vspace{2mm}

[KoZo] M. Kontsevich, A. Zorich, Connected components of the moduli spaces of Abelian differentials with prescribed singularities, arXiv:math/0201292.

\vspace{2mm}

[Moch] T. Mochizuki, Donaldson type invariants for algebraic surfaces. Lecture Notes in Mathematics vol. 1972, Springer, 2009.

\vspace{2mm}
[Nag1] K. Nagao, Derived categories of small toric Calabi-Yau 3-folds and counting invariants, arXiv:0809.2994.

\vspace{2mm}

[Nag2] K. Nagao,  Refined open noncommutative Donaldson-Thomas invariants for small crepant resolutions,  arXiv:0907.3784.
\vspace{2mm}

[NagNak] K. Nagao, H. Nakajima, Counting invariant of perverse coherent sheaves and its wall-crossing, arXiv:0809.2992.

\vspace{2mm}

[NiSe] J. Nicaise, J. Sebag, Motivic Serre invariants, ramification, and teh analytic Milnor fiber, ArXiv:math/0703217.

\vspace{2mm}

[OoYa] H. Ooguri, M. Yamazaki, Crystal Melting and Toric Calabi-Yau Manifolds, arXiv:0811.2801.

\vspace{2mm}
[Re1] M. Reineke, Cohomology of quiver moduli, functional equations, and integrality of Donaldson-Thomas type invariants, arXiv:0903.0261.

\vspace{2mm}

[So1] Y. Soibelman,  On non-commutative analytic spaces over non-archimedean fields,  arXiv:math/0606001.

\vspace{2mm}

[So2] Y. Soibelman, Two lectures on deformation theory and mirror symmetry, Journal of Math. Physics, vol.45:10, 2004.
\vspace{2mm}

[Sz] B. Szendroi,  Non-commutative Donaldson-Thomas theory and the conifold, arXiv:0705.3419.

[PT1] R. Pandharipande, R. P. Thomas,  Curve counting via stable pairs in the derived category,  arXiv:0707.2348.

\vspace{2mm}

[PT2] R. Pandharipande, R. P. Thomas, Stable pairs and BPS invariants, arXiv:0711.3899.
\vspace{2mm}

[Te1] M. Temkin,  Functorial desingularization of quasi-excellent schemes in characteristic zero: the non-embedded case, arXiv:0904.1592.
\vspace{2mm}

[Te2] M. Temkin,  Desingularization of quasi-excellent schemes in characteristic zero, arXiv:math/0703678.

\vspace{2mm}

[Te3] M. Temkin, Functorial desingularization over $\Q$: boundaries and the embedded case, arXiv:0912.2570.
\vspace{2mm}

[Zo] A.Zorich, Flat Surfaces, arXiv:math/0609392.

\vspace{3mm}
Addresses:

M.K.: IHES, 35 route de Chartres, F-91440, France

{maxim@ihes.fr}\\

Y.S.: Department of Mathematics, KSU, Manhattan, KS 66506, USA

{soibel@math.ksu.edu}

\end{document}